\documentstyle[amscd,amssymb,xypic,verbatim,11pt]{amsart}
\xyoption{all}

\newcommand{\nc}{\newcommand}
\newcommand{\rnc}{\renewcommand}

\nc{\Ups}{\Upsilon}
\nc{\Upss}{\Upsilon^{\mathsf S}}
\nc{\Upso}{\Upsilon^{\mathsf O}}

\nc{\exto}[1]{\stackrel{#1}{\longrightarrow}}
\nc{\dlim}{{\mathop{\lim\limits_{\longrightarrow}}}}
\nc{\lan}{\big\langle}
\nc{\ran}{\big\rangle}

\nc{\kk}{{\mathsf{k}}}
\nc{\ix}{{\mathsf{i}}}
\nc{\jx}{{\mathsf{j}}}

\nc{\C}{{\mathbb{C}}}
\nc{\HH}{{\mathbb{H}}}
\nc{\PP}{{\mathbb{P}}}
\nc{\QQ}{{\mathbb{Q}}}
\nc{\ZZ}{{\mathbb{Z}}}

\nc{\CA}{{\mathcal{A}}}
\nc{\CB}{{\mathcal{B}}}
\nc{\CC}{{\mathcal{C}}}
\nc{\D}{{\mathcal{D}}}
\nc{\CE}{{\mathcal{E}}}
\nc{\CF}{{\mathcal{F}}}
\nc{\CG}{{\mathcal{G}}}
\nc{\CH}{{\mathcal{H}}}
\nc{\CJ}{{\mathcal{J}}}
\nc{\CL}{{\mathcal{L}}}
\nc{\CM}{{\mathcal{M}}}
\nc{\CN}{{\mathcal{N}}}
\nc{\CO}{{\mathcal{O}}}
\nc{\CQ}{{\mathcal{Q}}}
\nc{\CR}{{\mathcal{R}}}
\nc{\CS}{{\mathcal{S}}}
\nc{\CT}{{\mathcal{T}}}
\nc{\CU}{{\mathcal{U}}}
\nc{\CV}{{\mathcal{V}}}
\nc{\CW}{{\mathcal{W}}}
\nc{\CX}{{\mathcal{X}}}
\nc{\CY}{{\mathcal{Y}}}
\nc{\CZ}{{\mathcal{Z}}}
\nc{\CMo}{{\mathcal{M}^\circ}}
\nc{\Co}{{{C}^\circ}}

\nc{\BY}{{\overline{Y}}}
\nc{\BYD}{{\overline{Y}{}^{|D|}}}
\nc{\OZ}{{\overline{Z}}}
\nc{\bg}{{\bar{g}}}

\nc{\bq}{{\mathbf{q}}}
\nc{\BD}{{\mathbf{D}}}
\nc{\BG}{{\mathbf{G}}}
\nc{\BM}{{\mathbf{M}}}
\nc{\BP}{{\mathbf{P}}}
\nc{\BZ}{{\mathbf{Z}}}
\nc{\BPr}{{\mathsf{P}}}
\nc{\BR}{{\mathbf{R}}}
\nc{\BRO}[1]{{{\mathbf{R}}^{\circ}_{#1}}}
\nc{\BRD}[1]{{{\mathbf{R}}^{|D|}_{#1}}}
\nc{\BRP}[1]{{{\mathbf{R}}^{1}_{#1}}}
\nc{\BRTP}[1]{{{\mathbf{\tilde{R}}}{}^{1}_{#1}}}
\nc{\BS}{{\mathbf{S}}}
\nc{\BMS}{{{\mathbf{M}}^{{s}}}}
\nc{\BMSS}{{{\mathbf{M}}^{{ss}}}}
\nc{\BMZ}{{\mathbf{M}^{\circ}}}
\nc{\BCL}{{\mathbf{L}}}

\nc{\PCC}{{{}^\perp\CC}}

\nc{\Cl}{{\mathsf{Cliff}}}
\nc{\Clev}{{\mathop{\mathsf{Cliff}}^{\circ}}}

\nc{\FA}{{\mathfrak{A}}}
\nc{\FB}{{\mathfrak{B}}}
\nc{\FF}{{\mathfrak{F}}}
\nc{\FI}{{\mathfrak{I}}}
\nc{\FZ}{{\mathfrak{Z}}}

\nc{\TFA}{{\tilde{\mathfrak{A}}}}
\nc{\TFB}{{\tilde{\mathfrak{B}}}}

\nc{\fa}{{\mathfrak{a}}}
\nc{\fg}{{\mathfrak{g}}}
\nc{\fp}{{\mathfrak{p}}}
\nc{\FD}{{\mathfrak{D}}}
\nc{\FE}{{\mathfrak{E}}}
\nc{\FL}{{\mathfrak{L}}}
\nc{\FM}{{\mathfrak{M}}}
\nc{\FR}{{\mathfrak{R}}}
\nc{\FS}{{\mathsf{S}}}

\nc{\sfc}{{\mathsf{c}}}
\nc{\sfch}{{\mathsf{ch}}}
\nc{\sfh}{{\mathsf{h}}}

\nc{\SK}{{\mathsf{K}}}
\nc{\SO}{{\mathsf{O}}}
\nc{\SQ}{{\mathsf{Q}}}
\nc{\SPV}{{\mathsf{S}^+\mathsf{V}}}
\nc{\SMV}{{\mathsf{S}^-\mathsf{V}}}
\nc{\SPMV}{{\mathsf{S}^\pm\mathsf{V}}}
\nc{\SX}{{S_X}}
\nc{\SY}{{S_Y}}
\nc{\phipsi}{{q}}
\nc{\eps}{\varepsilon}

\nc{\pim}{{\pi_-}}
\nc{\pip}{{\pi_+}}

\nc{\BE}{{\overline{\CE}}}
\nc{\TE}{{\tilde{\CE}}}
\nc{\TQ}{{\tilde{Q}}}
\nc{\TCF}{{\tilde{\CF}}}
\nc{\TCG}{{\tilde{\CG}}}
\nc{\TCL}{{\tilde{\CL}}}
\nc{\TF}{{\tilde{F}}}
\nc{\TW}{{\tilde{W}}}
\nc{\TCB}{{\widetilde{\CB}}}
\nc{\TCC}{{\tilde{\CC}}}
\nc{\TCX}{{\tilde{\CX}}}
\nc{\TCY}{{\tilde{\CY}}}
\nc{\TPhi}{{\tilde{\Phi}}}
\nc{\OPhi}{{\bar{\Phi}}}
\nc{\txi}{{\tilde{\xi}}}
\nc{\tp}{{\tilde{p}}}
\nc{\tq}{{\tilde{q}}}
\nc{\tzeta}{{\tilde{\zeta}}}
\nc{\tpi}{{\tilde{\pi}}}
\nc{\TUps}{{\tilde{\Ups}}}
\nc{\TUpss}{\TUps{}^{\mathsf S}}
\nc{\TUpso}{\TUps{}^{\mathsf O}}

\nc{\HCB}{{\widehat{\CB}}}
\nc{\HE}{{\widehat{\CE}}}
\nc{\HS}{{\widehat{S}}}
\nc{\HX}{{\hat{X}}}
\nc{\hxi}{{\hat{\xi}}}

\nc{\UH}{{\mathcal{H}}}

\nc{\TM}{{\widetilde{M}}}
\nc{\TCM}{{\widetilde{\CM}}}
\nc{\TS}{{\widetilde{S}}}
\nc{\TU}{{\widetilde{U}}}
\nc{\TX}{{\widetilde{X}}}
\nc{\TY}{{\widetilde{Y}}}
\nc{\TYO}{{{\widetilde{Y}}^\circ}}
\nc{\barf}{{\bar{f}}}
\nc{\te}{{\tilde{e}}{}}
\nc{\tf}{{\tilde{f}}}
\nc{\tg}{{\tilde{g}}}
\nc{\ti}{{\tilde{\imath}}}
\nc{\tj}{{\tilde{\jmath}}}
\nc{\ty}{{\tilde{y}}}
\nc{\tphi}{{\tilde{\phi}}}
\nc{\hf}{{\hat{f}}}

\nc{\urho}{{\underline{\rho}}}

\nc{\LRA}{\Leftrightarrow}
\nc{\RA}{\Rightarrow}
\nc{\lotimes}{\mathbin{\mathop{\otimes}\limits^{\mathbb{L}}}}
\nc{\CEnd}{\mathop{\mathcal{E}\mathit{nd}}\nolimits}
\nc{\CExt}{\mathop{\mathcal{E}\mathit{xt}}\nolimits}
\nc{\CHom}{\mathop{\mathcal{H}\mathit{om}}\nolimits}
\nc{\RH}{\mathop{{\mathsf{R}}\Gamma}\nolimits}
\nc{\RGamma}{\mathop{{\mathsf{R}}\Gamma}\nolimits}
\nc{\RHom}{\mathop{\mathsf{RHom}}\nolimits}
\nc{\RCHom}{\mathop{\mathsf{R}\mathcal{H}\mathit{om}}\nolimits}
\nc{\RG}{\mathop{\mathsf{R\Gamma}}\nolimits}
\nc{\Hom}{\mathop{\mathsf{Hom}}\nolimits}
\nc{\Ext}{\mathop{\mathsf{Ext}}\nolimits}
\nc{\End}{\mathop{\mathsf{End}}\nolimits}
\nc{\Tor}{\mathop{\mathsf{Tor}}\nolimits}
\nc{\Tordim}{\mathop{\mathsf{Tor}\text{\rm-}\mathsf{dim}}\nolimits}
\nc{\Hilb}{\mathop{\mathsf{Hilb}}\nolimits}
\nc{\Spec}{\mathop{\mathsf{Spec}}\nolimits}
\nc{\Proj}{\mathop{\mathsf{Proj}}\nolimits}
\nc{\Pic}{\mathop{\mathsf{Pic}}\nolimits}

\nc{\Tw}{\mathop{\mathsf{Tw}}\nolimits}
\nc{\Ker}{\mathop{\mathsf{Ker}}\nolimits}
\nc{\Coker}{\mathop{\mathsf{Coker}}\nolimits}
\nc{\codim}{\mathop{\mathsf{codim}}\nolimits}
\nc{\sing}{{\mathsf{sing}}}
\nc{\supp}{\mathop{\mathsf{supp}}}
\nc{\perf}{{\mathsf{perf}}}
\nc{\rank}{\mathop{\mathsf{rank}}}
\nc{\Pf}{{\mathsf{Pf}}}
\nc{\Gr}{{\mathsf{Gr}}}
\nc{\OGr}{{\mathsf{OGr}}}
\nc{\SGr}{{\mathsf{SGr}}}
\nc{\Flag}{{\mathsf{Fl}}}
\nc{\Kosz}{{\mathsf{Kosz}}}
\nc{\LGr}{{\mathsf{LGr}}}
\nc{\GTGr}{{\mathsf{G_2Gr}}}
\nc{\GTF}{{\mathsf{G_2F}}}
\nc{\OF}{{\mathsf{OF}}}
\nc{\Fl}{{\mathsf{Fl}}}
\nc{\Bl}{{\mathsf{Bl}}}
\nc{\GL}{{\mathsf{GL}}}
\nc{\PGL}{{\mathsf{PGL}}}
\nc{\SL}{{\mathsf{SL}}}
\nc{\SP}{{\mathsf{Sp}}}
\nc{\Spin}{{\mathsf{Spin}}}
\nc{\Tot}{{\mathsf{Tot}}}
\nc{\ev}{{\mathsf{ev}}}
\nc{\od}{{\mathsf{odd}}}
\nc{\coev}{{\mathsf{coev}}}
\nc{\id}{{\mathsf{id}}}
\nc{\opp}{{\mathsf{opp}}}
\nc{\PS}{{{\PP^3}}}
\nc{\Qu}{{{Q^3}}}
\nc{\tdim}{\mathop{\Tor\dim}}
\nc{\ecart}{{\fbox{$\scriptstyle\mathsf{EC}$}}}
\nc{\ad}{{\mathop{\mathsf ad}}}
\nc{\gr}{{\mathop{\mathsf gr}}}
\nc{\qgr}{{\mathop{\mathsf qgr}}}
\nc{\tor}{{\mathop{\mathsf tor}}}
\rnc{\mod}{{\mathop{\mathsf mod}}}
\nc{\Mod}{{\mathop{\mathsf Mod}}}
\nc{\Coh}{{\mathop{\mathsf Coh}}}
\nc{\Ab}{{\mathop{\mathcal{A}\mathit{b}}}}
\nc{\QCoh}{{\mathop{\mathsf QCoh}}}

\nc{\AAV}{{\mathcal{AAV}}}

\nc{\Rep}{{\mathsf{Rep}}}

\nc{\Cubics}{{{\mathcal{S}}_3}}
\nc{\VFT}{{{\mathcal{S}}_{14}}}
\nc{\VFTE}{{{\mathcal{N}}_{\mathrm{reg,sm}}}}
\nc{\MX}{{\CM_X}}
\nc{\MY}{{\CM_Y}}
\nc{\MYE}{{\CM_{Y,\CE}}}
\nc{\Yd}{{Y_d}}
\nc{\Yfive}{{Y_5}}
\nc{\Xg}{{X_{2g-2}}}
\nc{\Xtt}{{X_{22}}}
\nc{\Xst}{{X_{16}}}
\nc{\Xtw}{{X_{12}}}
\nc{\Xe}{{X_{8}}}
\nc{\Xf}{{X_{4}}}

\nc{\git}{{/\!\!/\!{}_\chi}}

\theoremstyle{plain}

\newtheorem{theorem}{Theorem}[section]
\newtheorem{conjecture}[theorem]{Conjecture}
\newtheorem{lemma}[theorem]{Lemma}
\newtheorem{proposition}[theorem]{Proposition}
\newtheorem{corollary}[theorem]{Corollary}

\theoremstyle{definition}

\newtheorem{definition}[theorem]{Definition}

\newtheorem{example}[theorem]{Example}

\theoremstyle{remark}

\newtheorem{remark}[theorem]{Remark}

\newenvironment{proof}{\noindent{\sf Proof:}}{\qed\medskip}

\title[Exceptional collections for Grassmannians of isotropic lines]%
{Exceptional collections for Grassmannians of isotropic lines}
\author{Alexander Kuznetsov}
\address{
Algebra Section, Steklov Mathematical Institute,
8 Gubkin str., Moscow 119991 Russia
}
\email{akuznet@@mi.ras.ru}
\date{}
\thanks{I was partially supported by RFFI grants 05-01-01034 and 02-01-01041,
Russian Presidential grant for young scientists No. MK-3926.2004.1,
CRDF Award No. RUM1-2661-MO-05, and the Russian Science Support Foundation.}

\begin{document}

\begin{abstract}
We construct a full exceptional collection of vector bundles in the
derived category of coherent sheaves on the Grassmannian of
isotropic two-dimensional subspaces in a symplectic vector space of
dimension $2n$ and in the derived category of coherent sheaves on
the Grassmannian of isotropic two-dimensional subspaces in an
orthogonal vector space of dimension $2n + 1$ for all~$n$.
\end{abstract}

\maketitle

\section{Introduction}

The derived category of coherent sheaves is the most important
algebraic invariant of an algebraic variety. This is a reason to
investigate its structure. In general, the structure of a
triangulated category (the derived category is triangulated) is
quite complicated. However, there is an important case when it can
be described fairly explicitly.

\begin{definition}[\cite{B,GR}]
A collection of objects $(E_1,E_2,\dots,E_n)$ of a $\kk$-linear triangulated category $\CT$ is {\sf exceptional}\/
if
$$
\text{$\RHom(E_i,E_i) = \kk$\quad for all $i$,\qquad and\qquad $\RHom(E_i,E_j) = 0$\quad for all $i>j$.}
$$
A collection $(E_1,E_2,\dots,E_n)$ is {\sf full}\/ if the minimal triangulated subcategory
of $\CT$ containing $E_1,E_2,\dots,E_n$ coincides with $\CT$.
\end{definition}

Triangulated categories possessing a full exceptional collection are the simplest among the others.
Every object of a triangulated category $\CT$ with a full exceptional collection
$(E_1,E_2,\dots,E_n)$ admits a unique (functorial) filtration with $i$-th quotient being a direct
sum of shifts of $E_i$. So, an exceptional collection can be considered as a kind of basis for triangulated category.
However, the condition of existence of a full exceptional collection in a triangulated category
is quite restrictive. For example, a necessary (but not sufficient) condition for the derived
category of coherent sheaves on a smooth projective variety $X$ to have a full exceptional collection
is the vanishing of nondiagonal Hodge numbers, that is $h^{ij}(X)=0$ for $i\ne j$
(an example of a smooth projective variety with vanishing nondiagonal Hodge numbers not admitting
a full exceptional collection is provided by Enriques surface).

The simplest example of a variety with a full exceptional collection is a projective space.
A.~Beilinson in~1978 showed~(\cite{Bei}) that the collection of line bundles
$(\CO,\CO(1),\dots,\CO(n))$ on $\PP^n$ is a full exceptional collection.
In~1988 M.~Kapranov constructed~(\cite{Ka}) full exceptional collections on Grassmannians
and flag varieties of groups ${\mathsf SL}_n$ and on smooth quadrics.
It has been conjectured afterwards that

\begin{conjecture}\label{c1}
Any projective homogeneous space of a semisimple algebraic group
admits a full exceptional collection consisting of vector bundles.
\end{conjecture}

\begin{conjecture}\label{c2}
The {\rm(}complete{\rm)} flag variety of a semisimple algebraic group
admits a full exceptional collection consisting of line bundles.
\end{conjecture}

It is somehow surprising that only a very little progress in this direction has been achieved.
As it was already mentioned Kapranov showed that both~\ref{c1} and~\ref{c2} are true
for the group~${\mathsf SL}_n$. Also it is easy to see that~\ref{c2} is true for the group ${\mathsf SP}_{2n}$
(\cite{S2}) since the corresponding flag variety can be represented as an iterated projectivization
of vector bundles, and by~\cite{O} the derived category of a projectivization of a vector bundle admits
a full exceptional collection consisting of line bundles if the base does.

If one wishes to establish conjecture~\ref{c1}, i.e. to construct a full exceptional collection on any homogeneous space,
it is natural to consider first the case of Grassmannians, i.e. of the homogeneous varieties $G/P$, where $P$ is a maximal
parabolic subgroup in a semisimple algebraic group $G$. Presumably, this would suffice to prove conjecture~\ref{c1}
in general via the parabolic induction procedure. Now let us briefly describe what is known about
derived categories of Grassmannians of classical semisimple algebraic groups.

{\bf (${\mathsf SL}_n$)} The Grassmannians of the group ${\mathsf SL}_n$ are the usual Grassmannians $\Gr(k,n)$ of $k$-dimensional
subspaces in an $n$-dimensional vector space for $1\le k\le n-1$. Let $\CU$ denote the tautological
rank $k$ subbundle in the trivial vector bundle $\CO^{\oplus n}_{\Gr(k,n)}$.
Consider the corresponding principal $\GL_k$-bundle on $\Gr(k,n)$.
Given a nonincreasing collection of $k$ integers $\alpha = (\alpha_1\ge\alpha_2\ge\dots\ge\alpha_k)$
we consider it as a dominant weight of the group ${\mathsf GL}_k$
and denote by $\Sigma^\alpha\CU$ the associated vector bundle on $\Gr(k,n)$. In particular,
\begin{itemize}
\item[] for $\alpha = (m,0,\dots,0)$, $\Sigma^\alpha\CU = S^m\CU$ is the $m$-th symmetric power, and
\item[] for $\alpha = (\underbrace{1,\dots,1}_{m},0,\dots,0)$, $\Sigma^\alpha\CU = \Lambda^m\CU$
is the $m$-th exterior power.
\end{itemize}
Kapranov has shown \cite{Ka} that
\begin{itemize}
\item for all $1\le k\le n-1$ the collection of vector bundles
$
\{ \Sigma^\alpha\CU^*\ |\ n-k \ge \alpha_1\ge\dots\ge\alpha_k\ge 0 \}
$
is a full exceptional collection on $\Gr(k,n)$.
\end{itemize}

{\bf (${\mathsf SP}_{2n}$)} The Grassmannians of the group ${\mathsf SP}_{2n}$ are
the isotropic Grassmannians of $k$-dimensional isotropic subspaces in a symplectic
$2n$-dimensional vector space for $1\le k\le n$ which we denote by $\SGr(k,2n)$.
For these Grassmannians exceptional collections are known only in the following cases:
\begin{itemize}
\item $k=1$: in this case $\SGr(1,2n) \cong \PP^{2n-1}$ and we have the Beilinson's exceptional collection;
\item $n=3$: exceptional collections on $\SGr(2,6)$ and $\SGr(3,6)$ were constructed by A.Samokhin~\cite{S2,S1}.
\end{itemize}
In the present paper we construct a full exceptional collection on $\SGr(2,2n)$ for all~$n$.

{\bf (${\mathsf SO}_n$)} The Grassmannians of the group ${\mathsf SO}_{n}$ are
the isotropic Grassmannians of $k$-dimensional isotropic subspaces in an orthogonal
$n$-dimensional vector space for $1\le k < n/2$ if $n=2m+1$ is odd, and $1\le k\le n/2$, $k\ne n/2-1$,
if $n = 2m$ is even which we denote by $\OGr(k,n)$ (if $n=2m$ then $\OGr(m,2m)$ has two connected
components, each of them is a Grassmannian and $\OGr(m-1,2m)$ is a projectivization
of a vector bundle over a component of $\OGr(m,2m)$, so it is not a Grasmannian).
For these Grassmannians exceptional collections are known only in the following cases:
\begin{itemize}
\item $k=1$ and odd $n$: in this case $\OGr(1,n) \cong Q^{n-2}$ is an odd-dimensional quadric and\\
$\{\CO_Q,S,\CO_Q(1),\CO_Q(2),\dots,\CO_Q(n-3)\}$ is a full exceptional collection~\cite{Ka} ($S$ is the spinor bundle);
\item $k=1$ and even $n$: in this case $\OGr(1,n) \cong Q^{n-2}$ is an even-dimensional quadric and\\
$\{\CO_Q,S^+,S^-,\CO_Q(1),\CO_Q(2),\dots,\CO_Q(n-3)\}$ is a full exceptional collection~\cite{Ka} ($S^\pm$ are the spinor bundles);
\item $n \le 6$: the group ${\mathsf SO}_n$ in this case has a simpler description:
${\mathsf SO}_3$ up to a finite subgroup coincides with ${\mathsf SL}_2$,
${\mathsf SO}_4$ up to a finite subgroup coincides with ${\mathsf SL}_2 \times {\mathsf SL}_2$,
${\mathsf SO}_5$ up to a finite subgroup coincides with ${\mathsf SP}_4$,
and
${\mathsf SO}_6$ up to a finite subgroup coincides with ${\mathsf SL}_4$.
\end{itemize}
In the present paper we construct a full exceptional collection on $\OGr(2,2n+1)$ for all~$n$.

As we already mentioned above the main result of the present paper is a construction
of full exceptional collections on Grassmannians $\SGr(2,2m)$ and $\OGr(2,2m+1)$.
Let us say some words about the method of constructing and proving fullness of these collections.
First of all we consider a special exceptional collection on $\Gr(2,2m)$:
\begin{equation}\label{lcgr2m}
\left(
\begin{array}{rrrrrrr}
S^{m-1}\CU^* & S^{m-1}\CU^*(1) & \dots & S^{m-1}\CU^*(m-1) \\
S^{m-2}\CU^* & S^{m-2}\CU^*(1) & \dots & S^{m-2}\CU^*(m-1) & S^{m-2}\CU^*(m) & \dots & S^{m-2}\CU^*(2m-1) \\
\vdots\quad & \vdots\quad & & \vdots\qquad\qquad & \vdots\qquad & & \vdots\qquad\qquad \\
S^{2}\CU^* & S^{2}\CU^*(1) & \dots & S^{2}\CU^*(m-1) & S^{2}\CU^*(m) & \dots & S^{2}\CU^*(2m-1) \\
\CU^* & \CU^*(1) & \dots & \CU^*(m-1) & \CU^*(m) & \dots & \CU^*(2m-1) \\
\CO & \CO(1) & \dots & \CO(m-1) & \CO(m) & \dots & \CO(2m-1)
\end{array}\right)
\end{equation}
(the ordering is from bottom to top in columns, the left column goes first).
Note that this is {\em not}\/ the Kapranov's collection. In comparison
with the Kapranov's collection, \eqref{lcgr2m} is more symmetric with respect
to the $\CO(1)$-twisting. This symmetry property is axiomatized in the notion
of a {\em minimal Lefschetz exceptional collection}.

A Lefschetz exceptional collection is just an exceptional collection which
consists of several blocks, each of them is a subblock of the previous one
twisted by $\CO(1)$ (in the collection~\eqref{lcgr2m} the blocks are the columns).
A Lefschetz exceptional collection is minimal, if, roughly speaking,
its first block is minimal possible (see the precise definition in~\cite{K2}).
The Kapranov's collection on $\Gr(2,2m)$ is a Lefschetz collection with the first block
$(\CO,\CU^*,S^2\CU^*,\dots,S^{2m-2}\CU^*)$, however it is not minimal, since the first
block of the collection~\eqref{lcgr2m} is strictly smaller.

Lefschetz exceptional collections enjoy a lot of nice properties,
for example they behave well with respect to the restriction to a hyperplane section.
Explicitly, if we remove the first block of a Lefschetz exceptional collection
and then restrict the rest to a hyperplane, we obtain a Lefschetz exceptional
collection (see Proposition~\ref{lhp}). Evidently, the smaller is the first block,
the bigger collection on the hyperplane we obtain. In particular, the biggest exceptional
collection on a hyperplane can be obtained from a minimal Lefschetz collection
on the ambient variety.

Since $\SGr(2,2m)$ is a hyperplane
section of $\Gr(2,2m)$ we obtain in this way a Lefschetz exceptional collection for $\SGr(2,2m)$
\begin{equation}\label{lcsgr}
\left(
\begin{array}{rrrrrr}
S^{m-1}\CU^* & \dots & S^{m-1}\CU^*(m-2) \\
S^{m-2}\CU^* & \dots & S^{m-2}\CU^*(m-2) & S^{m-2}\CU^*(m-1) & \dots & S^{m-2}\CU^*(2m-2) \\
\vdots\quad & & \vdots\qquad\qquad & \vdots\qquad & & \vdots\qquad\qquad \\
S^{2}\CU^* & \dots & S^{2}\CU^*(m-2) & S^{2}\CU^*(m-1) & \dots & S^{2}\CU^*(2m-2) \\
\CU^* & \dots & \CU^*(m-2) & \CU^*(m-1) & \dots & \CU^*(2m-2) \\
\CO & \dots & \CO(m-2) & \CO(m-1) & \dots & \CO(2m-2)
\end{array}\right)
\end{equation}
(we applied to the obtained collection on $\SGr(2,2m)$ an additional $\CO(-1)$-twist for convenience).
This collection turns out to be full.
To prove this we use the fact that this exceptional collection behaves well
with respect to the restriction to any $\SGr(2,2m-2) \subset \SGr(2,2m)$, and
use induction in $m$.

Similarly, in the case of the orthogonal Grassmannian $\OGr(2,2m+1)$ we start with a Lefschetz
exceptional collection on $\Gr(2,2m+1)$:
\begin{equation}\label{lcgr2mp1}
\left(
\begin{array}{rrrrrrr}
S^{m-1}\CU^* & S^{m-1}\CU^*(1) & \dots & S^{m-1}\CU^*(m-1) & S^{m-1}\CU^*(m) & \dots & S^{m-1}\CU^*(2m) \\
S^{m-2}\CU^* & S^{m-2}\CU^*(1) & \dots & S^{m-2}\CU^*(m-1) & S^{m-2}\CU^*(m) & \dots & S^{m-2}\CU^*(2m) \\
\vdots\quad & \vdots\quad & & \vdots\qquad\qquad & \vdots\qquad & & \vdots\qquad\qquad \\
S^{2}\CU^* & S^{2}\CU^*(1) & \dots & S^{2}\CU^*(m-1) & S^{2}\CU^*(m) & \dots & S^{2}\CU^*(2m) \\
\CU^* & \CU^*(1) & \dots & \CU^*(m-1) & \CU^*(m) & \dots & \CU^*(2m) \\
\CO & \CO(1) & \dots & \CO(m-1) & \CO(m) & \dots & \CO(2m)
\end{array}\right)
\end{equation}
This time $\OGr(2,2m+1) \subset \Gr(2,2m+1)$ is not a hyperplane section,
rather it is a section of a vector bundle $S^2\CU^*$. Nevertheless,
a subcollection of~\eqref{lcgr2mp1} restricts to a Lefschetz exceptional
collection on $\OGr(2,2m+1)$. However, in a contrast with the symplectic case
it is not full. Similarly, as in the case of quadrics we need suitably defined
spinor bundles on $\OGr(2,2m+1)$. These bundles can be defined as follows.

Consider the partial flag variety $\OF(2,m;2m+1) \subset \OGr(2,2m+1) \times \OGr(m,2m+1)$
and denote by $\pi_2:\OF(2,m;2m+1) \to \OGr(2,2m+1)$, $\pi_m:\OF(2,m;2m+1) \to \OGr(m,2m+1)$
the projections. We define the spinor bundle $\CS$ on $\OGr(2,2m+1)$ as
$$
\CS = \pi_{2*}\pi_m^*\CO_{\OGr(m,2m+1)}(1),
$$
where $\CO_{\OGr(m,2m+1)}(1)$ is the positive generator of $\Pic(\OGr(m,2m+1))$
(which is {\em not }\/ the pullback of $\CO_{\Gr(2,2m+1)}(1)$,
actually the pullback of $\CO_{\Gr(2,2m+1)}(1)$ is isomorphic to
$\CO_{\OGr(m,2m+1)}(2)$). It turns out that the twists of spinor
bundles can be inserted in the exceptional collection, so that the collection
\begin{equation}\label{lcogr2mp1}
\left(
\begin{array}{rrrrrrr}
\CS & \CS(1) & \dots & \CS(m-1) & \CS(m) & \dots & \CS(2m-3) \\
S^{m-2}\CU^* & S^{m-2}\CU^*(1) & \dots & S^{m-2}\CU^*(m-1) & S^{m-2}\CU^*(m) & \dots & S^{m-2}\CU^*(2m-3) \\
\vdots\quad & \vdots\quad & & \vdots\qquad\qquad & \vdots\qquad & & \vdots\qquad\qquad \\
S^{2}\CU^* & S^{2}\CU^*(1) & \dots & S^{2}\CU^*(m-1) & S^{2}\CU^*(m) & \dots & S^{2}\CU^*(2m-3) \\
\CU^* & \CU^*(1) & \dots & \CU^*(m-1) & \CU^*(m) & \dots & \CU^*(2m-3) \\
\CO & \CO(1) & \dots & \CO(m-1) & \CO(m) & \dots & \CO(2m-3)
\end{array}\right)
\end{equation}
is exceptional.
To prove its fullness we use the fact that this exceptional collection behaves well
with respect to the restriction to any $\OGr(2,2m-1) \subset \OGr(2,2m+1)$, and
use induction in $m$.

\medskip

It also worth mentioning that the notion of a Lefschetz exceptional collection
(or more generally of a Lefschetz semiorthogonal decomposition) was introduced in~\cite{K1}
as a starting point for the theory of Homological Projective Duality. Therefore it is natural to ask
what will be the Homologically Projectively Dual for $\Gr(2,n)$, $\SGr(2,2m)$ and $\OGr(2,2m+1)$
with respect to the considered here Lefschetz exceptional collections.
In~\cite{K3} we answer this question for the Grassmannians $\Gr(2,6)$ and $\Gr(2,7)$ and formulate
a conjectural answer for $\Gr(2,n)$ for all $n$. Since the symplectic Grassmannian $\SGr(2,2m)$ is a hyperplane
section of $\Gr(2,2m)$, \cite{K3} also gives an answer for $\SGr(2,2m)$.
For the orthogonal Grassmannian $\OGr(2,2m+1)$ the question is still open.

The paper is organized as follows.
In section~\ref{s_lc} we remind the definition and some properties of Lefschetz exceptional collections and
in section~\ref{s_bbw} we remind the Borel-Bott-Weil Theorem.
In section~\ref{s_gr} we show that~\eqref{lcgr2m} and~\eqref{lcgr2mp1} are full exceptional collections
on Grassmannians $\Gr(2,2m)$ and $\Gr(2,2m+1)$ respectively,
and give a proof of fullness of these collections by induction in $m$,
demonstrating a method used later in the case of isotropic Grassmannians.
In section~\ref{s_sgr} we show that~\eqref{lcsgr} is a full exceptional collection
on the symplectic isotropic Grassmannians $\SGr(2,2m)$.
In section~\ref{s_sb} we develop a theory of spinor bundles on homogeneous spaces
of orthogonal groups.
Finally, in section~\ref{s_ogr2mp1} we show that~\eqref{lcogr2mp1} is a full exceptional collection
on the orthogonal isotropic Grassmannians $\OGr(2,2m+1)$.

\section{Lefschetz exceptional collections}\label{s_lc}

Let $X$ be a smooth projective algebraic variety over a field $\kk$ with an ample line bundle $\CO_X(1)$.
We denote by $\D^b(X)$ the bounded derived category of coherent sheaves on $X$.


\begin{definition}[cf.~\cite{K1}]
A {\sf Lefschetz collection}\/ in $\D^b(X)$ with respect to the line bundle $\CO_X(1)$
is a collection of objects of $\D^b(X)$ which has a block structure
$$
\Big(
\underbrace{E_1,E_2,\dots,E_{\lambda_0}}_{\text{block $1$}},
\underbrace{E_1(1),E_2(1),\dots,E_{\lambda_1}(1)}_{\text{block $2$}},
\dots,
\underbrace{E_1(\ix-1),E_2(\ix-1),\dots,E_{\lambda_{\ix-1}}(\ix-1)}_{\text{block $\ix$}}
\Big),
$$
where $\lambda = (\lambda_0 \ge \lambda_1 \ge \dots \ge \lambda_{\ix-1} > 0)$ is a nonincreasing
sequence of positive integers (the {\sf support partition} of the Lefschetz collection).
\end{definition}
In other words, a collection is Lefschetz with support partition $\lambda  = (\lambda_0,\lambda_1,\dots,\lambda_{\ix-1})$
if it splits into $\ix$ blocks of length $\lambda_0$, $\lambda_1$, \dots, $\lambda_{\ix-1}$ such that the $k$-th block
consists of the first $\lambda_{k-1}$ objects of the first block twisted by $\CO_X(k-1)$.
A Lefschetz decomposition is uniquely determined by its first block $(E_1,E_2,\dots,E_{\lambda_0})$ and
its support partition $\lambda$.


Recall that a collection $(E_1,E_2,\dots,E_n)$ of objects in $\D^b(X)$ is called {\sf exceptional} if
\begin{enumerate}
\item for each $i$ the object $E_i$ is exceptional, i.e. $\Hom(E_i,E_i) = \kk$,
$\Ext^p(E_i,E_i) = 0$ for $p\ne 0$;
\item for all $i < j$ we have $\Ext^\bullet(E_j,E_i) = 0$.
\end{enumerate}

The symmetry of Lefschetz collections simplifies the verification of exceptionality.

\begin{lemma}\label{critl}
A Lefschetz collection $E_\bullet$ with support partition $\lambda  = (\lambda_0,\lambda_1,\dots,\lambda_{\ix-1})$ is exceptional
iff
\begin{enumerate}
\item its first block $(E_1,E_2,\dots,E_{\lambda_0})$ is an exceptional collection, and
\item $\Ext^\bullet(E_p,E_q(-k)) = 0$ {for $1\le k\le \ix-1$, $1\le p\le \lambda_k$, and $1\le q\le \lambda_0$}.
\end{enumerate}
\end{lemma}
\begin{proof}
It suffices to note that $\Ext^\bullet(E_p(k),E_q(l)) = \Ext^\bullet(E_p,E_q(l-k))$.
\end{proof}

\begin{example}
Now let us give several examples of Lefschetz exceptional collections
\begin{itemize}
\item Any exceptional collection is a 1-block Lefschetz collection.
\item For any $d > 0$ the standard exceptional collection $(\CO_{\PP^n},\CO_{\PP^n}(1),\dots,\CO_{\PP^n}(n))$ on $\PP^n$
is Lefschetz with respect to $\CO_{\PP^n}(d)$ with support partition
$\lambda = (\underbrace{d,d,\dots,d}_q,r)$, where $n+1 = qd + r$, $0 < r\le d$.
\item The Kapranov's exceptional collection on an odd-dimensional quadric
$$
(\CO_{Q^n},\CS,\CO_{Q^n}(1),\dots,\CO_{Q^n}(n-1))
$$
is Lefschetz with respect to $\CO_{Q^n}(1)$
with support partition $\lambda = (2,\underbrace{1,\dots,1}_{n-1})$. Here $\CS$ is the spinor bundle.
\item The Kapranov's exceptional collection on an even-dimensional quadric
$$
(\CO_{Q^n},\CS^+,\CS^-,\CO_{Q^n}(1),\CO_{Q^n}(2),\dots,\CO_{Q^n}(n-1))
$$
is Lefschetz with respect to $\CO_{Q^n}(1)$ with support partition $\lambda = (3,\underbrace{1,\dots,1}_{n-1})$.
Here $\CS^+$ and $\CS^-$ are the spinor bundles on $Q^n$. However, this Lefschetz collection is not minimal.
Indeed, it is easy to see that the mutation of $\CS^-$ through $\CO_{Q^n}(1)$
is isomorphic to $\CS^+(1)$, hence the following collection on $Q^n$
$$
(\CO_{Q^n},\CS^+,\CO_{Q^n}(1),\CS^+(1),\CO_{Q^n}(2),\dots,\CO_{Q^n}(n-1))
$$
is also a full exceptional Lefschetz collection with respect to $\CO_{Q^n}(1)$
but with another support partition $\lambda' = (2,2,\underbrace{1,\dots,1}_{n-2})$.
\end{itemize}
\end{example}

A simple but very useful property is that a Lefschetz exceptional collection in $\D^b(X)$
with respect to $\CO_X(1)$ gives a Lefschetz exceptional collection
on any hyperplane section of $X$.
Explicitly, if we remove the first block of a Lefschetz exceptional collection
and then restrict the rest to a hyperplane, we obtain a Lefschetz exceptional
collection collection.

\begin{proposition}[\cite{K1}]\label{lhp}
Let
$(E_1,\dots,E_{\lambda_0},E_1(1),\dots,E_{\lambda_1}(1),\dots,E_1(\ix-1),\dots,E_{\lambda_{\ix-1}}(\ix-1))$
be a Lefschetz exceptional collection in $\D^b(X)$ with support partition $\lambda$.
Let $Y \subset X$ be a hyperplane section with respect to $\CO_X(1)$.
Then
$(E_1,\dots,E_{\lambda_1},E_1(1),\dots,E_{\lambda_2}(1),\dots,E_1(\ix-2),\dots,E_{\lambda_{\ix-1}}(\ix-2))$
is a Lefschetz exceptional collection in $\D^b(Y)$ with support partition $\lambda' = (\lambda_1,\lambda_2,\dots,\lambda_{\ix-1})$.
\end{proposition}
\begin{proof}
We have $\Ext_Y^\bullet(E_p,E_q(-k)) = H^\bullet(Y,E_p^*\otimes E_q(-k))$.
Since $Y$ is a hyperplane section of $X$ with respect to $\CO_X(1)$ we have a resolution
$$
0 \to \CO_X(-1) \to \CO_X \to \CO_Y \to 0.
$$
Tensoring it with $E_p^*\otimes E_q(-k)$ we obtain a long exact sequence
$$
\dots \to
H^s(X,E_p^*\otimes E_q(-k)) \to
H^s(Y,E_p^*\otimes E_q(-k)) \to
H^{s+1}(X,E_p^*\otimes E_q(-k-1)) \to
\dots
$$
Taking $k=0$ we deduce that $\Ext^\bullet_Y(E_p,E_q) = \Ext^\bullet_X(E_p,E_q)$
for $1 \le p,q \le \lambda_1$, hence the collection $(E_1,E_2,\dots,E_{\lambda_1})$ on $Y$ is exceptional.
Taking $1\le k\le \ix-2$ we deduce that $\Ext^\bullet(E_p,E_q(-k)) = 0$ for
$1\le p\le \lambda_{k+1}$ and $1\le q\le \lambda_1$.
By lemma~\ref{critl} this means that the desired Lefschetz collection on $Y$ is exceptional.
\end{proof}

\section{The Borel--Bott--Weil Theorem}\label{s_bbw}

The Borel--Bott--Weil Theorem computes the cohomology of line bundles
on the flag variety of a semisimple algebraic group. It also can be used to compute
the cohomology of equivariant vector bundles on Grassmannians.
We restrict here to the case of the group $\GL(V)$
(see, however, remark~\ref{anygroup}).

Let $V$ be a vector space of dimension $n$. The standard identification
of the weight lattice of the group $\GL(V)$ with $\ZZ^n$ takes the $k$-th
fundamental weight $\pi_k$ (the heighest weight of the representation $\Lambda^kV$)
to the vector $(1,1,\dots,1,0,0,\dots,0) \in \ZZ^n$ (the first $k$ entries
are 1, and the last $n-k$ are 0). Under this identification
the cone of dominant weights of $\GL(V)$ gets identified with the set of nonincreasing
sequences $\alpha = (a_1,a_2,\dots,a_n)$ of integers. For such $\alpha$
we denote by $\Sigma^\alpha V = \Sigma^{a_1,a_2,\dots,a_n} V$ the corresponding
representation of $\GL(V)$.
Note that $\Sigma^{1,1,\dots,1}V = \det V$.

Similarly, given a vector bundle $E$ of rank $n$ on a scheme $S$
we consider the corresponding principal $\GL(n)$-bundle on $S$ and
denote by $\Sigma^\alpha E$ the vector bundle associated
with the $\GL(n)$-representation of highest weight $\alpha$.

The group $\BS_n$ of permutations acts naturally on the weight lattice $\ZZ^n$.
Denote by $\ell:\BS_n \to \ZZ$ the standard length function.
Note that for every $\alpha \in \ZZ^n$ there exists a permutation $\sigma \in \BS_n$
such that $\sigma(\alpha)$ is nonincreasing. If all entries of $\alpha$
are distinct then such $\sigma$ is unique and $\sigma(\alpha)$ is strictly decreasing.

Let $X$ be the flag variety of $\GL(V)$. Let $L_\alpha$ denote the line bundle
on $X$ corresponding to the weight~$\alpha$ (so that $L_{\pi_k}$ is the pullback
of $\CO_{\PP(\Lambda^kV)}(1)$ under the natural projection $X \to \PP(\Lambda^kV)$).

Denote by
$$
\rho = (n,n-1,\dots,2,1)
$$
half the sum of the positive roots of $\GL(V)$.
The corresponding line bundle $L_\rho$ is the square root of the anticanonical line bundle.

The Borell-Bott-Weil Theorem computes the cohomology of line bundles $L_\alpha$ on $X$.
\begin{theorem}[\cite{D}]\label{bbw}
Assume that all entries of $\alpha + \rho$ are distinct.
Let $\sigma$ be the unique permutation such that $\sigma(\alpha+\rho)$ is strictly decreasing.
Then
$$
H^k(X,L_\alpha) = \begin{cases}
\Sigma^{\sigma(\alpha+\rho) - \rho}V^*, & \text{if $k = \ell(\sigma)$}\\
0, & \text{otherwise}
\end{cases}
$$
If not all entries of $\alpha + \rho$ are distinct then $H^\bullet(X,L_\alpha) = 0$.
\end{theorem}

\begin{remark}\label{anygroup}
The Borel--Bott-Weil Theorem is true for any semisimple algebraic group.
One should replace in the statement $\ZZ^n$ by the weight lattice,
the set of strictly decreasing sequences by the interior of the dominant cone,
$\rho$ by half the sum of the positive roots, and the group $\BS_n$ by the Weil group.
\end{remark}

Now consider a Grassmannian $\Gr(k,V)$. Let $\CU \subset V\otimes\CO_{\Gr(k,V)}$ denote
the tautological subbundle of rank $k$. Denote by $W/\CU$ the corresponding quotient bundle
and by $\CU^\perp$ its dual, so that we have the following (mutually dual) exact sequences
$$
0 \to \CU \to V\otimes\CO_{\Gr(k,V)} \to V/\CU \to 0,
\qquad
0 \to \CU^\perp \to V^*\otimes\CO_{\Gr(k,V)} \to \CU^* \to 0.
$$
Note that $\Sigma^{1,1,\dots,1}\CU^* \cong \Sigma^{-1,-1,\dots,-1}\CU^\perp$ is the positive
generator of $\Pic{\Gr(k,V)}$.
Let $\pi:X \to {\Gr(k,V)}$ denote the canonical projection from the flag variety to the Grassmannian.

\begin{proposition}[\cite{Ka}]\label{equb}
Let $\beta \in \ZZ^k$ and $\gamma \in \ZZ^{n-k}$ be nonincreasing sequences.
Let $\alpha = (\beta,\gamma) \in \ZZ^n$ be their concatenation.
Then we have $\pi_* L_\alpha \cong \Sigma^\beta\CU^*\otimes\Sigma^\gamma\CU^\perp$.
\end{proposition}

\begin{corollary}\label{bbwg}
If $\beta \in \ZZ^k$ and $\gamma \in \ZZ^{n-k}$ are nonincreasing sequences
and $\alpha = (\beta,\gamma) \in \ZZ^n$ then
$$
H^\bullet({\Gr(k,V)},\Sigma^\beta\CU^*\otimes\Sigma^\gamma\CU^\perp) \cong H^\bullet(X,L_\alpha).
$$
\end{corollary}

Since every $\GL(V)$-equivariant vector bundle on ${\Gr(k,V)}$ is isomorphic to
$\Sigma^\beta\CU^*\otimes\Sigma^\gamma\CU^\perp$ for some nonincreasing $\beta\in\ZZ^k$,
$\gamma\in\ZZ^{n-k}$, a combination of corollary~\ref{bbwg} with theorem~\ref{bbw}
allows to compute the cohomology of any equivariant vector bundle on ${\Gr(k,V)}$.

In a combination with the Littlewood--Richardson rules the Borel--Bott--Weil Theorem
allows to compute $\Ext$ groups between equivariant bundles on the Grassmannians.
As an example we will prove the following lemma which will be used later in section~\ref{s_gr}.

\begin{lemma}\label{bbwgr2n}
Let $X = \Gr(2,V)$, $\dim V = n$. If $0 \le l_1,l_2 \le n/2 - 1$ and $0\le k\le n-1$ then we have
$$
\Ext^p(S^{l_1}\CU^*,S^{l_2}\CU^*(-k)) =
\begin{cases}
S^{l_2-l_1}V^*, & \text{if $l_1 \le l_2$, $k = 0$, and $p = 0$}\\
\kk, & \text{if $l_1 = l_2 = n/2 - 1$, $k = n/2$, and $p = n - 2$}\\
0, & \text{otherwise}
\end{cases}
$$
\end{lemma}
\begin{proof}
First of all, we have
$$
\Ext^\bullet(S^{l_1}\CU^*,S^{l_2}\CU^*(-k)) =
H^\bullet(X,S^{l_1}\CU\otimes S^{l_2}\CU^*(-k)).
$$
Further, by the Littlewood-Richardson rule we have
$$
S^{l_1}\CU\otimes S^{l_2}\CU^*(-k) =
\bigoplus_{t=0}^{\text{\sf min}\{l_1,l_2\}} \Sigma^{l_2-k-t,-l_1-k+t}\CU^* =
\Sigma^{l_2-k,-l_1-k}\CU^* \oplus S^{l_1-1}\CU\otimes S^{l_2-1}\CU^*(-k),
$$
so it suffices to compute $H^\bullet(X,\Sigma^\alpha\CU^*)$ for $\alpha = (l_2-k,-l_1-k,0,0,\dots,0)$.
Note that $$
\alpha + \rho = (n+l_2-k,n-l_1-k-1,n-2,n-3,\dots,1).
$$
Assume that all entries of $\alpha + \rho$ are distinct. This is equivalent to
$$
n+l_2-k \not\in \{n-2,n-3,\dots,1\}
\qquad\text{and}\qquad
n-l_1-k-1 \not\in \{n-2,n-3,\dots,1\}.
$$
The second condition implies that either $n-l_1-k-1 = n-1$, that is $l_1 = k = 0$,
or $n-l_1-k-1\le 0$, that is $l_1 + k \ge n - 1$.
In the first case we get $\alpha + \rho = (n+l_2,n-1,n-2,\dots,1)$,
$\sigma = 1$, hence $H^0(X,\Sigma^\alpha\CU^*) = S^{l_2}V^*$ and other cohomologies are zero.
In the second case, since $l_1 \le n/2 - 1$ it follows that $k \ge n/2$.
Since $l_2 \le n/2 - 1$ we have $n+l_2-k \le n - 1$,
and since $k\le n-1$ we have $n+l_2-k \ge 1$.
Therefore, the first of the above conditions implies that
$n+l_2-k = n-1$, hence $l_1 = l_2 = n/2 - 1$ and $k = n/2$.
In this case $\alpha + \rho = (n-1,0,n-2,n-3,\dots,1)$,
$\ell(\sigma) = n - 2$, $\sigma(\alpha+\rho)-\rho = (-1,-1,\dots,-1)$,
hence $H^{n-2}(X,\Sigma^\alpha\CU^*) = \kk$ and other cohomologies are zero.
Combining all this we deduce the lemma.
\end{proof}

\section{Usual Grassmannian}\label{s_gr}

Consider the Grassmannian $X = \Gr(2,W)$ of two-dimensional subspaces in an $n$-dimensional
vector space $W$. Let $\CU$ denote the tautological rank $2$ subbundle on $X = \Gr(2,W)$.
We will distinguish between the cases of even and odd $n$. Let
$$
m = \left\lfloor \frac{n}2 \right\rfloor,
$$
so that either $n=2m$ or $n=2m+1$.
If $n = 2m$ we consider the collection~\eqref{lcgr2m} on $X$ and
if $n = 2m + 1$ we consider the collection~~\eqref{lcgr2mp1} on $X$.
These are Lefschetz collections with the first block
\begin{equation}\label{ugc}
(\CO_X,\CU^*,S^2\CU^*,\dots,S^{m-1}\CU^*)
\end{equation}
and with the support partition
\begin{equation}\label{ugp}
\lambda = \begin{cases}
(\underbrace{m,m,\dots,m}_{2m+1}), & \text{if $n = 2m+1$}\\
(\underbrace{m,m,\dots,m}_{m},\underbrace{m-1,m-1,\dots,m-1}_{m}), & \text{if $n = 2m$}
\end{cases}
\end{equation}
In the other words, the collections consist of vector bundles $S^l\CU^*(k)$ with integers $(k,l)$
from the set
\begin{equation}\label{Sn}
\Ups_n = \{(k,l)\in\ZZ^2\ |\ \text{$0\le k\le n-1$, $0\le l\le m-1$ and $l\le m-2$ for $k\ge m$ and even $n$}\}.
\end{equation}
It is interesting to compare these collections with the standard Kapranov's collections. The latter also consist
of vector bundles $S^l\CU^*(k)$ but with other restrictions on possible values of $(k,l)$, namely
$0\le k,l$ and $k+l\le n-2$. On the following picture we draw the triangles corresponding to Kapranov's
exceptional collections together with the regions $\Ups_n$.
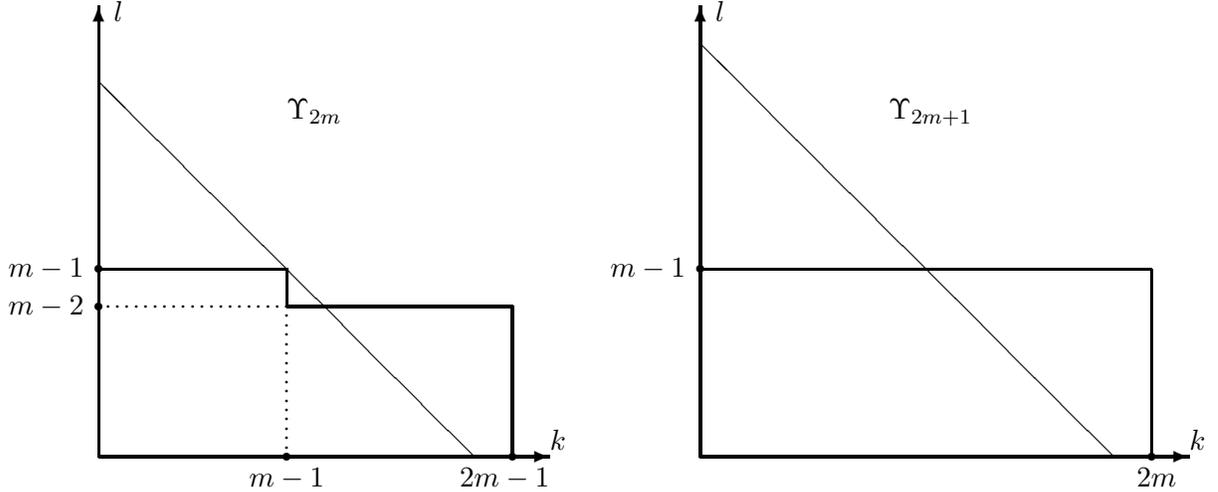
\begin{figure}[h]
\unitlength=1mm
\begin{picture}(170,65)
\thicklines
\put(45,50){$\Ups_{2m}$}
\put(20,5){\vector(0,1){60}}
\put(20,5){\vector(1,0){60}}
\put(20,30){\line(1,0){25}}
\put(45,25){\line(0,1){5}}
\put(45,25){\line(1,0){30}}
\put(75,5){\line(0,1){20}}
\multiput(45,25)(-1.25,0){21}{\circle*{.1}}
\multiput(45,25)(0,-1.25){17}{\circle*{.1}}
\put(20,30){\circle*{1}}
\put(8,29){$m-1$}
\put(20,25){\circle*{1}}
\put(8,24){$m-2$}
\put(45,5){\circle*{1}}
\put(40,1){$m-1$}
\put(75,5){\circle*{1}}
\put(68,1){$2m-1$}
\put(80,6){$k$}
\put(22,63){$l$}
%
%
\put(125,50){$\Ups_{2m+1}$}
\put(100,5){\vector(0,1){60}}
\put(100,5){\vector(1,0){65}}
\put(100,30){\line(1,0){60}}
\put(160,5){\line(0,1){25}}
\put(100,30){\circle*{1}}
\put(88,29){$m-1$}
\put(160,5){\circle*{1}}
\put(158,1){$2m$}
\thinlines
\put(70,5){\line(-1,1){50}}
\put(155,5){\line(-1,1){55}}
\put(165,6){$k$}
\put(102,63){$l$}
\end{picture}
\end{figure}

The main result of this section is the following

\begin{theorem}
Let $X = \Gr(2,n)$ and let $\Ups_n$ be the set defined in~\eqref{Sn}.
Then the Lefschetz collection $\{S^l\CU^*(k)\ |\ (k,l)\in\Ups_n\}$
is a full exceptional collection in $\D^b(X)$.
\end{theorem}

Certainly, the collection $\{S^l\CU^*(k)\ |\ (k,l)\in\Ups_n\}$ can be obtained
by a sequence of mutations from the Kapranov's exceptional collection
(this is really not too complicated, one should consider complexes~\eqref{crucial}
constructed below). However, we prefer to use an inductive argument, since it can
(and will) be applied for the symplectic and orthogonal Grassmannians as well.

For a start we must check that the collection is exceptional.
This is easily done by the Borel--Bott--Weil theorem.

\begin{lemma}\label{upsexc}
The Lefschetz collection $\{S^l\CU^*(k)\ |\ (k,l)\in\Ups_n\}$ is exceptional in $\D^b(X)$.
\end{lemma}
\begin{proof}
Combine lemma~\ref{bbwgr2n} and lemma~\ref{critl}.
%
%
\end{proof}

It remains to prove the fullness of the collection.
We start with some preparations.

Besides the set $\Ups_n$ consider also the set
\begin{equation}\label{TSn}
\hspace{-10pt}
\TUps_{n-1} = \{(k,l)\in\ZZ^2\ |\ \text{$0\le k\le n-1$, $0\le l\le m-1$ and $l\le m-2$ for $k\ge m-1$ and odd $n$}\}.
\end{equation}
On the following picture the region $\Ups_n$ is drawn together with the region $\TUps_{n-1}$.
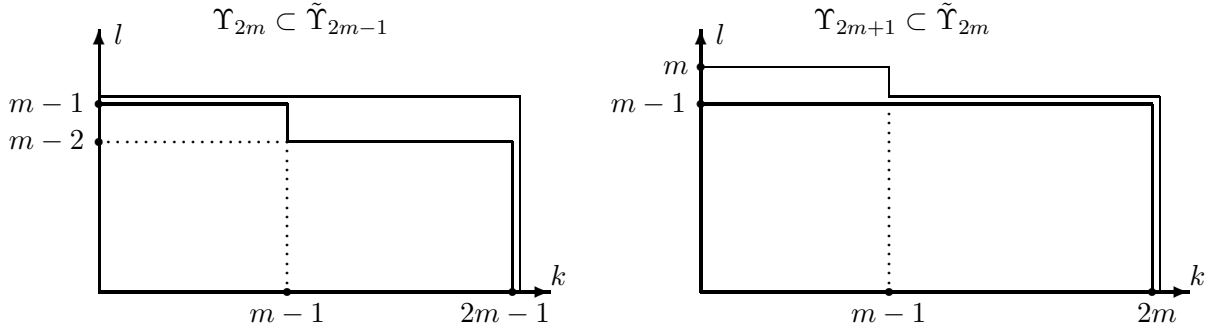
\begin{figure}[h]
\unitlength=1mm
\begin{picture}(170,45)
\thicklines
\put(35,40){$\Ups_{2m} \subset\TUps_{2m-1}$}
\put(20,5){\vector(0,1){35}}
\put(20,5){\vector(1,0){60}}
\put(20,30){\line(1,0){25}}
\put(45,25){\line(0,1){5}}
\put(45,25){\line(1,0){30}}
\put(75,5){\line(0,1){20}}
\multiput(45,25)(-1.25,0){21}{\circle*{.1}}
\multiput(45,25)(0,-1.25){17}{\circle*{.1}}
\put(20,30){\circle*{1}}
\put(8,29){$m-1$}
\put(20,25){\circle*{1}}
\put(8,24){$m-2$}
\put(45,5){\circle*{1}}
\put(40,1){$m-1$}
\put(75,5){\circle*{1}}
\put(68,1){$2m-1$}
\put(80,6){$k$}
\put(22,38){$l$}
%
%
\put(115,40){$\Ups_{2m+1} \subset\TUps_{2m}$}
\put(100,5){\vector(0,1){35}}
\put(100,5){\vector(1,0){65}}
\put(100,30){\line(1,0){60}}
\put(160,5){\line(0,1){25}}
\multiput(125,30)(0,-1.25){21}{\circle*{.1}}
\put(100,30){\circle*{1}}
\put(88,29){$m-1$}
\put(125,5){\circle*{1}}
\put(120,1){$m-1$}
\put(100,35){\circle*{1}}
\put(95,34){$m$}
\put(160,5){\circle*{1}}
\put(158,1){$2m$}
\put(165,6){$k$}
\put(102,38){$l$}
\thinlines
\put(20,31){\line(1,0){56}}
\put(76,5){\line(0,1){26}}
\put(100,35){\line(1,0){25}}
\put(125,31){\line(0,1){4}}
\put(125,31){\line(1,0){36}}
\put(161,5){\line(0,1){26}}

\end{picture}
\end{figure}

\begin{lemma}\label{tsigma}
For any $(k,l)\in\TUps_{n-1}$ the vector bundle $S^l\CU^*(k)$ lies in the triangulated category
generated by the Lefschetz collection $\{S^l\CU^*(k)\ |\ (k,l)\in\Ups_n\}$ in $\D^b(\Gr(2,W))$.
%
\end{lemma}
\begin{proof}
Let $n=2m$. Then $\TUps_{n-1} \setminus \Ups_n = \{(m,m-1),(m+1,m-1),\dots,(2m-1,m-1)\}$.
So, we have to check that the vector bundles $S^{m-1}\CU^*(m)$, $S^{m-1}\CU^*(m+1)$, \dots, $S^{m-1}\CU^*(2m-1)$
lie in the triangulated category generated by the Lefschetz collection $\{S^l\CU^*(k)\ |\ (k,l)\in\Ups_n\}$ on $\Gr(2,W)$.

Consider the (dual) tautological exact sequence $0 \to\CU^\perp \to W^*\otimes\CO_X \to \CU^* \to 0$
on $X = \Gr(2,W)$. It induces the following long exact sequence
\begin{equation}\label{skus}
0 \to \Lambda^k \CU^\perp \to \Lambda^k W^*\otimes\CO_X \to
\Lambda^{k-1}W^*\otimes \CU^* \to \dots \to W^*\otimes S^{k-1}\CU^* \to
S^k\CU^* \to 0
\end{equation}
for any $0 \le k \le n-2$.
Dualizing, using isomorphisms $\CU^{\perp*} \cong W/\CU$, $S^l\CU
\cong S^l\CU^*\otimes\CO_X(-l)$ and replacing $k$ by $n-2-k$ we
obtain exact sequence
\begin{multline}\label{sku}
0 \to S^{n-2-k}\CU^*(k+2-n) \to W\otimes S^{n-3-k}\CU^*(k+3-n) \to
\dots \to \\
\to \Lambda^{n-3-k}W\otimes \CU^*(-1) \to
\Lambda^{n-2-k}W\otimes\CO_X \to \Lambda^{n-2-k}(W/\CU) \to 0.
\end{multline}
On the other hand, we have an isomorphism
$$
\Lambda^k\CU^\perp \cong \Lambda^{n-2-k}(W/\CU)\otimes\CO_X(-1).
$$
Using this isomorphism for gluing the sequence~(\ref{skus}) twisted
by $n-k-1$ with the sequence~(\ref{sku}) twisted by $n-k-2$ we
obtain the following exact sequence
\begin{multline}\label{crucial}
0 \to S^{n-2-k}\CU^* \to
\\ \to
W\otimes S^{n-3-k}\CU^*(1) \to \dots \to
\Lambda^{n-3-k}W\otimes \CU^*(n-k-3) \to
\Lambda^{n-2-k}W\otimes\CO_X(n-k-2) \to
\\ \to
\Lambda^kW^*\otimes\CO_X(n-k-1) \to \Lambda^{k-1}W^*\otimes
\CU^*(n-k-1) \to \dots \to W^*\otimes S^{k-1}\CU^*(n-k-1) \to
\\ \to
S^k\CU^*(n-k-1) \to 0.
\end{multline}

Take $k=m-1$. Then the sequence~(\ref{crucial}) gives a decomposition
for $S^{m-1}\CU^*(m)$ with respect to the Lefschetz collection.
Twisting this sequence by $\CO_X(1),\dots,\CO_X(m-1)$ we obtain also decompositions
for $S^{m-1}\CU^*(m+1),\dots,S^{m-1}\CU^*(2m-1)$.

Now let $n=2m+1$. Then $\TUps_{n-1} \setminus \Ups_n = \{(0,m),(1,m),\dots,(m-1,m)\}$.
Take $k=m-1$. Then the sequence~(\ref{crucial}) gives a decomposition
for $S^m\CU^*$ with respect to the Lefschetz collection.
Twisting this sequence by $\CO_X(1),\dots,\CO_X(m-1)$ we obtain also decompositions
for $S^m\CU^*(1),\dots,S^m\CU^*(m-1)$.
\end{proof}

Another preparatory result is the following:

\begin{lemma}\label{xphi}
For any $0 \ne \phi \in W^* = H^0(\Gr(2,W),\CU^*)$ the zero locus of $\phi$ on $X = \Gr(2,W)$
is the Grassmannian $X_\phi = \Gr(2,\Ker\phi) \subset \Gr(2,W) = X$. Moreover, we have
the following resolution of the structure sheaf $\CO_{X_\phi}$ on $X$:
\begin{equation}\label{resxphi}
0 \to \CO_X(-1) \to \CU \to \CO_X \to i_{\phi*}\CO_{X_\phi} \to 0,
\end{equation}
where $i_{\phi*}:X_\phi \to X$ is the embedding.
\end{lemma}
\begin{proof}
The first part is evident. For the second part we note that any nonzero
section $\phi$ of $\CU^*$ is regular since $\dim X_\phi = 2(n-1)-4 = \dim X - 2$,
so the sheaf $i_{\phi*}\CO_{X_\phi}$ admits a Koszul resolution which takes form~(\ref{resxphi}).
\end{proof}

Now we are ready for the proof of the theorem.
We use induction in $n$. The base of induction, $n=3$, is clear.
Indeed, in this case $X = \Gr(2,W) = \PP^2$ and the Lefschetz collection
takes form $(\CO_X,\CO_X(1),\CO_X(2))$ which is well known to be full.

Now assume that the fullness of the corresponding Lefschetz collection is already proved for $n-1$.
Assume also that the Lefschetz collection for $n$ is not full. Then by~\cite{B} there exists an object $F\in\D^b(X)$,
right orthogonal to all bundles in the collection, hence by lemma~\ref{tsigma} to all $S^l\CU^*(k)$
with $(k,l) \in \TUps_{n-1}$, i.e.
$$
0 = \RHom(S^l\CU^*(k),F) = H^\bullet(X,S^l\CU(-k)\otimes F)
\qquad\text{for all $(k,l)\in\TUps_{n-1}$}
$$
Let us check that $i_\phi^*F = 0$ for any $0 \ne \phi\in W^*$. For this we take any $(k,l)\in \Ups_{n-1}$
and tensor the resolution~(\ref{resxphi}) by $S^l\CU(-k)\otimes F$. Taking into account isomorphism
$$
(S^l\CU(-k)\otimes F)\otimes i_{\phi*}\CO_{X_\phi} \cong
i_{\phi*}i_\phi^*(S^l\CU(-k)\otimes F) \cong
i_{\phi*}(S^l\CU(-k)\otimes i_\phi^*(F))
$$
we get a resolution
$$
0 \to S^l\CU(-k-1)\otimes F \to S^l\CU(-k)\otimes \CU \otimes F \to S^l\CU(-k)\otimes F \to i_{\phi*}(S^l\CU(-k)\otimes i_\phi^*(F))) \to 0.
$$
Now note, that for $(k,l)\in\Ups_{n-1}$ we have
$$
(k+1,l),(k,l+1),(k,l) \in \TUps_{n-1},
\qquad\text{and also $(k+1,l-1)\in\TUps_{n-1}$ if $l\ge 1$}.
$$
Since $S^l\CU(-k) \otimes\CU = S^{l+1}\CU(-k) \oplus S^{l-1}\CU(1-k)$ (the second summand vanishes if $l=0$),
it follows that the cohomology on $X$ of the first three terms of the above complex vanishes.
Therefore we have
$$
\RHom_{X_\phi}(S^l\CU^*(k),i_\phi^*F) = H^\bullet(X_\phi,S^l\CU(-k)\otimes i_\phi^*(F)) = 0
\qquad\text{for all $(k,l)\in\Ups_{n-1}$}
$$
Thus $i_\phi^* F$ lies in the right orthogonal to the subcategory of $\D^b(X_\phi)$
generated by the exceptional collection $\{S^l\CU^*(k)\ |\ (k,l)\in\Ups_{n-1}\}$
which by the induction hypothesis is full. Hence indeed $i_\phi^*F = 0$.
So we conclude by the following

\begin{lemma}
If for $F\in\D^b(X)$ we have $i_\phi^*F = 0$ for any $0\ne\phi\in W^*$ then $F = 0$.
\end{lemma}
\begin{proof}
Assume that $F \ne 0$. Let $q$ be the maximal integer such that $\CH^q(F)\ne 0$,
take a point $x\in\supp\CH^q(F)$ and choose $0\ne \phi\in W^*$ such that $x\in X_\phi$
(this is equivalent to the vanishing of a linear function $\phi$ on the 2-dimensional subspace
of $W$ corresponding to $x\in X = \Gr(2,W)$). Since the functor $i_\phi^*$ is left-exact
it easily follows that $\CH^q(i_\phi^* F) \ne 0$, so $i_\phi^*F \ne 0$.
\end{proof}

Thus we have proved that the desired collection is indeed full.

\section{Symplectic Grassmannian}\label{s_sgr}

Consider the isotropic Grassmannian $X = \SGr(2,W)$ of two-dimensional subspaces in
a symplectic vector space $W$ of dimension $2m$.
Note that $\SGr(2,W)$ is a hyperplane section of the usual Grassmannian $\Gr(2,W)$,
the hyperplane coressponds to the symplectic form $\omega \in \Lambda^2W^* = H^0(\Gr(2,W),\CO(1))$.
Let $\CU$ denote the restriction of the tautological rank $2$ subbundle from $\Gr(2,W)$ to $X = \SGr(2,W)$.
Restricting the Lefschetz exceptional collection on $\Gr(2,W)$ constructed in the previous section
and using proposition~\ref{lhp} we obtain the Lefschetz collection (\ref{lcsgr}) on $\SGr(2,W)$.
This is a Lefschetz collection with the first block~(\ref{ugc}) and the support partition
\begin{equation}\label{ugps}
\lambda^{\mathsf S} = (\underbrace{m,m,\dots,m}_{m-1},\underbrace{m-1,m-1,\dots,m-1}_{m}).
\end{equation}
In the other words, the collection consists of vector bundles $S^l\CU^*(k)$ with integers $(k,l)$
from the set
\begin{equation}\label{Sns}
\Upss_{2m} = \{(k,l)\in\ZZ^2\ |\ \text{$0\le k\le 2m-2$, $0\le l\le m-1$ and $l\le m-2$ for $k\ge m-1$}\}.
\end{equation}
On the following picture we draw the regions $\Ups_{2m}$ and $\Upss_{2m}$
in the $(k,l)$-plane.
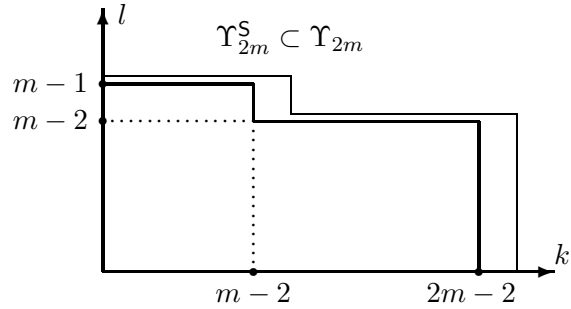
\begin{figure}[h]
\unitlength=1mm
\begin{picture}(170,45)
\thicklines
\put(75,35){$\Upss_{2m} \subset \Ups_{2m}$}
\put(60,5){\vector(0,1){35}}
\put(60,5){\vector(1,0){60}}
\put(60,30){\line(1,0){20}}
\put(80,25){\line(0,1){5}}
\put(80,25){\line(1,0){30}}
\put(110,5){\line(0,1){20}}
\put(60,30){\circle*{1}}
\put(48,29){$m-1$}
\put(60,25){\circle*{1}}
\put(48,24){$m-2$}
\put(80,5){\circle*{1}}
\put(75,1){$m-2$}
\put(110,5){\circle*{1}}
\put(103,1){$2m-2$}
\multiput(80,25)(-1.25,0){17}{\circle*{.1}}
\multiput(80,25)(0,-1.25){17}{\circle*{.1}}
\put(120,6){$k$}
\put(62,38){$l$}
\thinlines
\put(60,31){\line(1,0){25}}
\put(85,26){\line(0,1){5}}
\put(85,26){\line(1,0){30}}
\put(115,5){\line(0,1){21}}
%

\end{picture}
\end{figure}


\begin{theorem}
Let $\Upss_{2m}$ be the set defined in~\eqref{Sns}.
Then the Lefschetz collection $\{S^l\CU^*(k)\ |\ (k,l)\in\Upss_{2m}\}$
is a full exceptional collection in $\D^b(\SGr(2,2m))$.
\end{theorem}

As it was already mentioned above the collection is exceptional by proposition~\ref{lhp}.
So, it remains to prove the fullness of the collection.
We use the method demonstrated in the previous section.
The main difference is that the induction step now changes $\dim W$ by $2$.

Consider the set
\begin{equation}\label{TSns}
\TUpss_{2m-2} = \{(k,l)\in\ZZ^2\ |\ \text{$0\le k\le 2m-2$, $0\le l\le m$ and $l\le m-1$ for $k\ge m-2$}\}.
\end{equation}
On the following picture the region $\Upss_{2m}$ is drawn together with the region $\TUpss_{2m-2}$.

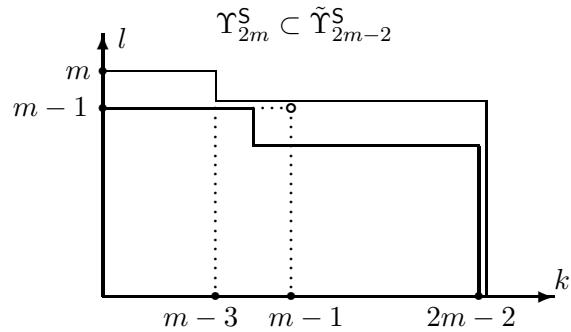
\begin{figure}[h]
\unitlength=1mm
\begin{picture}(170,45)
\thicklines
\put(75,40){$\Upss_{2m} \subset \TUpss_{2m-2}$}
\put(60,5){\vector(0,1){35}}
\put(60,5){\vector(1,0){60}}
\put(60,30){\line(1,0){20}}
\put(80,25){\line(0,1){5}}
\put(80,25){\line(1,0){30}}
\put(110,5){\line(0,1){20}}

\put(85,30){\circle{1}}

\thinlines
\put(60,35){\line(1,0){15}}
\put(75,31){\line(0,1){4}}
\put(75,31){\line(1,0){36}}
\put(111,5){\line(0,1){26}}

\put(60,30){\circle*{1}}
\put(48,29){$m-1$}
\put(60,35){\circle*{1}}
\put(55,34){$m$}
\put(75,5){\circle*{1}}
\put(68,1){$m-3$}
\put(85,5){\circle*{1}}
\put(82,1){$m-1$}
\put(110,5){\circle*{1}}
\put(103,1){$2m-2$}
\multiput(83.75,30)(-1.25,0){4}{\circle*{.1}}
\multiput(85,28.75)(0,-1.25){20}{\circle*{.1}}
\multiput(75,30)(0,-1.25){21}{\circle*{.1}}
\put(120,6){$k$}
\put(62,38){$l$}

\end{picture}
\end{figure}

\noindent
The small circle corresponds to the bundle $S^{m-1}\CU^*(m-1)$ which plays special role as it will be seen further.

\begin{lemma}\label{tsigmas}
For any $(k,l)\in\TUpss_{2m-2}$ the vector bundle $S^l\CU^*(k)$ lies in the triangulated category
generated by the Lefschetz collection $\{S^l\CU^*(k)\ |\ (k,l)\in\Upss_{2m}\}$ in $\D^b(\SGr(2,W))$.
%
\end{lemma}
\begin{proof}
Note that
$$
\TUpss_{2m-2} \setminus \Upss_{2m} = \{(0,m),(1,m),\dots,(m-3,m)\} \cup \{(m-1,m-1)\} \cup \{(m,m-1),\dots,(2m-2,m-1)\}.
$$
So, we have to check that the vector bundles
$$
S^m\CU^*,S^m\CU^*(1),\dots,S^m\CU^*(m-3),
\qquad
S^{m-1}\CU^*(m-1),
\qquad
S^{m-1}\CU^*(m),\dots,S^{m-1}\CU^*(2m-2)
$$
lie in the triangulated category generated by the Lefschetz collection $\{S^l\CU^*(k)\ |\ (k,l)\in\Upss_{2m}\}$.
Consider the restriction to $\SGr(2,W)\subset\Gr(2,W)$ of the exact sequence~(\ref{crucial}) with $k=m-2$.
It gives a decomposition of $S^m\CU^*$ with respect to the Lefschetz collection.
Twisting this sequence by $\CO_X(1),\dots,\CO_X(m-3)$ we obtain also decompositions
for $S^m\CU^*(1),\dots,S^m\CU^*(m-3)$. Similarly, taking $k=m-1$ we obtain a decomposition
of $S^{m-1}\CU^*(m)$ with respect to the Lefschetz collection. Twisting this sequence by
$\CO_X(1),\dots,\CO_X(m-2)$ we obtain also decompositions for $S^{m-1}\CU^*(m+1),\dots,S^{m-1}\CU^*(2m-2)$.
It remains only to find a decomposition for the vector bundle $S^{m-1}\CU^*(m-1)$.
This is done in the following proposition.
\end{proof}

\begin{proposition}
On $\SGr(2,W)$ there exists a bicomplex
$$
\scriptsize
\arraycolsep=.5pt
\begin{array}{ccccccccccccccccccccc}
S^{m-1}\CU^*
\\
\downarrow &&
\\
W^*\otimes S^{m-2}\CU^*(1) & \to &
S^{m-1}\CU^*(1)
\\
\downarrow &&
\downarrow &&
\\
\Lambda^2W^*\otimes S^{m-3}\CU^*(2) & \to &
W^*\otimes S^{m-2}\CU^*(2) & \to &
S^{m-1}\CU^*(2)
\\
\downarrow &&
\downarrow &&
\downarrow &&
%
\\
\vdots &&
\vdots &&
\vdots &&
\ddots &&
%
\\
\downarrow &&
\downarrow &&
\downarrow &&
&&
%
\\
\Lambda^{m-2}W^*\otimes\CU^*_X(m-2) & \to &
\Lambda^{m-3}W^*\otimes S^2\CU^*(m-2) & \to &
\Lambda^{m-4}W^*\otimes S^3\CU^*(m-2) & \to &
\dots & \to &
S^{m-1}\CU^*(m-2)
\\
\downarrow &&
\downarrow &&
\downarrow &&
&&
\downarrow &&
\\
\Lambda^{m-1}W^*\otimes\CO_X(m-1) & \to &
\Lambda^{m-2}W^*\otimes\CU^*(m-1) & \to &
\Lambda^{m-3}W^*\otimes S^2\CU^*(m-1) & \to &
\dots & \to &
W^*\otimes S^{m-2}\CU^*(m-1) & \to &
S^{m-1}\CU^*(m-1)
\end{array}
$$
the total complex of which is exact.
\end{proposition}
\begin{proof}
First, consider the above diagram on the ambient Grassmannian $\Gr(2,W)$
with rows being the twisted truncations of the exact sequences~(\ref{skus})
and columns being the twisted truncations of the exact sequences~(\ref{sku})
where we identify $\Lambda^kW$ with $\Lambda^kW^*$ via the symplectic form $\omega$.
Certainly, this diagram is not commutative, but let us check that it commutes modulo $\omega$.
More precisely, we will check that the compositions of arrows in the square
$$
\xymatrix{
\Lambda^{k}W^*\otimes S^{m-1-k}\CU^* \ar[r] \ar[d] & \Lambda^{k-1}W^*\otimes S^{m-k}\CU^* \ar[d] \\
\Lambda^{k+1}W^*\otimes S^{m-2-k}\CU^*(1) \ar[r]   & \Lambda^{k}W^*\otimes S^{m-1-k}\CU^* (1)
}
$$
which is a typical (up to a twist) square of the diagram, coincide modulo $\omega$.
First of all, we apply the functor $\Hom(S^{m-1-k}\CU^*,-)$ to this square
and check that the resulted square of $\Hom$-s commutes modulo~$\omega$.
Indeed, since by the Borel--Bott--Weil Theorem we have
$$
\begin{array}{ll}
\Hom(S^{m-1-k}\CU^*,S^{m-1-k}\CU^*) = \kk,
\qquad &
\Hom(S^{m-1-k}\CU^*,S^{m-k}\CU^*) = W^*,\\
\Hom(S^{m-1-k}\CU^*,S^{m-2-k}\CU^*(1)) = W^*,
\qquad &
\Hom(S^{m-1-k}\CU^*,S^{m-1-k}\CU^*(1)) = W^*\otimes W^*,
\end{array}
$$
the square of $\Hom$-s takes form
$$
\xymatrix{
\Lambda^{k}W^* \ar[rr]^-c \ar[d]^-{\omega c\omega^{-1}} &&
\Lambda^{k-1}W^*\otimes W^* \ar[d]^-{\omega c\omega^{-1}} \\
\Lambda^{k+1}W^*\otimes W^* \ar[rr]^-c &&
\Lambda^{k}W^*\otimes W^*\otimes W^*
}
$$
where $c$ is the canonical map and $\omega c\omega^{-1}$ is the canonical map conjugated by $\omega$.
Let $\{e_i\}$ be a base of $W$, and $\{f_i\}$ be the dual base of $W^*$.
Then the compositions of arrows in this square act as follows
$$
%
\xymatrix@R=5pt{
\alpha \ar@{|->}[rr]^-c &&
\sum (\alpha\vdash e_i)\otimes f_i \ar@{|->}[r]^-{\omega c\omega^{-1}} &
\sum \;(\alpha\vdash e_i)\wedge f_j \otimes  \omega(e_j) \otimes f_i,\\
\alpha \ar@{|->}[rr]^-{\omega c\omega^{-1}} &&
\sum (\alpha\wedge f_j)  \otimes  \omega(e_j) \ar@{|->}[r]^-c &
\sum [(\alpha\vdash e_i)\wedge f_j \otimes  \omega(e_j) \otimes f_i + (-1)^k \alpha\wedge (f_j\vdash e_i) \otimes f_i\otimes\omega(e_j)],
}
$$
where $\alpha \in \Lambda^kW^*$ and $\vdash$ denotes the convolution of a form and a vector.
It remains to note that we have $f_j\vdash e_i = \delta_{ij}$ and $\sum f_i\otimes\omega(e_i) = \omega$,
hence indeed the difference of the compositions is given by the map
$\Lambda^kW^* \to \Lambda^kW^*\otimes \Lambda^2W^* \to \Lambda^{k}W^*\otimes W^*\otimes W^*$,
$\alpha \mapsto \alpha\otimes\omega$.
Finally we note that the canonical maps
$\Hom(S^{m-1-k}\CU^*,F) \otimes S^{m-1-k}\CU^* \to F$
for $F = S^{m-1-k}\CU^*$, $S^{m-k}\CU^*$, $S^{m-2-k}\CU^*(1)$, and $S^{m-1-k}\CU^*(1)$
are surjective, hence it follows that the squares of the constructed diagram commute modulo $\omega$.

Now we restrict the diagram to the isotropic Grassmannian $\SGr(2,W)$.
It follows that the squares now commute (because $\omega$ vanishes),
so we got a bycomplex. It remains to check that its total complex is exact.

Consider the spectral sequences of this bycomplex.
The first term of the first spectral sequence (the cohomology of rows)
is concentrated at the left column by~(\ref{skus}), which means that
the total complex can have nontrivial cohomology only in the first $m$ terms.
On the other hand,
the first term of the second spectral sequence (the cohomology of columns)
is concentrated at the bottom row by~(\ref{sku}), which means that
the total complex can have nontrivial cohomology only in the last $m$ terms.
Combining these two observations we deduce that the total complex can have
nontrivial cohomology only in the middle ($m$-th) term.
Finally, considering again the second spectral sequence and using~(\ref{sku})
we see that this cohomology is a subsheaf of the sheaf $\Lambda^{m-1}(W/\CU)\otimes\CO_X(m-1)$,
hence is torsion free. But computing the Euler characteristics of the total complex
we see that its rank is zero, hence the cohomology vanishes.
\end{proof}

Another preparatory result is the following:

\begin{lemma}\label{xphis}
For any two-dimensional subspace $\langle w_1,w_2\rangle \subset W \cong W^* = H^0(\SGr(2,W),\CU^*)$
such that $\omega(w_1,w_2)\ne 0$, the zero locus of the corresponding section
$\phi = \phi_{w_1,w_2} \in H^0(\SGr(2,W),\CU^*\oplus\CU^*)$ on $X = \SGr(2,W)$
is the isotropic Grassmannian $X_{w_1,w_2} = \SGr(2,\lan w_1,w_2\ran^\perp) \subset \SGr(2,W) = X$.
Moreover, we have the following resolution of the structure sheaf $\CO_{X_\phi}$ on $X$:
\begin{equation}\label{resxphis}
\hspace{-5mm}
0 \to \CO_X(-2) \to \CU(-1)\oplus \CU(-1) \to \CO_X(-1)^{\oplus 3} \oplus S^2\CU \to
\CU\oplus\CU \to \CO_X \to i_{\phi*}\CO_{X_{w_1,w_2}} \to 0,
\end{equation}
where $i_{\phi}:X_{w_1,w_2} \to X$ is the embedding.
\end{lemma}
\begin{proof}
The first part is evident (one should only note that the restriction of $\omega$
to the subspace $\langle w_1,w_2\rangle^\perp$ in nondegenerate provided $\omega(w_1,w_2)\ne 0$).
For the second part we note that any such section $\phi = \phi_{w_1,w_2}$ of $\CU^*\oplus\CU^*$
is regular since $\dim X_{w_1,w_2} = 2(n-2)-5 = \dim X - 4$, so the sheaf $i_{\phi*}\CO_{X_{w_1,w_2}}$
admits a Koszul resolution which takes form~(\ref{resxphis}).
\end{proof}

Now we are ready for the proof of the theorem.
We use induction in $m$. The base of induction, $m=2$, is clear.
Indeed, in this case $X = \SGr(2,W) = Q^3$, a three-dimensional quadric,
and the Lefschetz collection takes form $(\CO_X,\CU^*,\CO_X(1),\CO_X(2))$
which is well known to be full (actually, this is precisely the Kapranov's
exceptional collection for $Q^3$).

Now assume that the fullness of the corresponding Lefschetz collection is already proved for $m-1$.
Assume also that the Lefschetz collection for $m$ is not full. Then by~\cite{B} there exists an object $F\in\D^b(X)$,
right orthogonal to all bundles in the collection, hence by lemma~\ref{tsigmas}
$$
0 = \RHom(S^l\CU^*(k),F) = H^\bullet(X,S^l\CU(-k)\otimes F)
\qquad\text{for all $(k,l)\in\TUpss_{2m-2}$}
$$
by lemma~\ref{tsigmas}.
Let us check that $i_\phi^*F = 0$ for any $\phi = \phi_{w_1,w_2}$ like in lemma~\ref{xphis}.
For this we take any $(k,l)\in \Upss_{2m-2}$ and tensor the resolution~(\ref{resxphis})
by $S^l\CU(-k)\otimes F$. Taking into account isomorphism
$$
(S^l\CU(-k)\otimes F)\otimes i_{\phi*}\CO_{X_\phi} \cong
i_{\phi*}(i_\phi^*(S^l\CU(-k)\otimes F)) \cong
i_{\phi*}(S^l\CU(-k)\otimes i_\phi^*(F))
$$
and noting that for $(k,l)\in\Upss_{2m-2}$ we have
$$
(k+2,l),(k+2,l-1),(k+1,l+1),(k+1,l),(k+2,l-2),(k,l+2),(k,l+1),(k+1,l-1),(k,l) \in \TUpss_{2m-2},
$$
it follows that the cohomology on $X$ of the first five terms of the above complex vanishes.
Therefore we have
$$
\RHom_{X_\phi}(S^l\CU^*(k),i_\phi^*F) = H^\bullet(X_\phi,S^l\CU(-k)\otimes i_\phi^*(F)) = 0
\qquad\text{for all $(k,l)\in\Upss_{2m-2}$}
$$
Thus $i_\phi^* F$ lies in the right orthogonal to the subcategory of $\D^b(X_\phi)$
generated by the exceptional collection $\{S^l\CU^*(k)\ |\ (k,l)\in\Upss_{2m-2}\}$
which by the induction hypothesis is full. Hence indeed $i_\phi^*F = 0$.
So we conclude by the following

\begin{lemma}
If for $F\in\D^b(X)$ we have $i_\phi^*F = 0$ for any two-dimensional subspace
$\langle w_1,w_2\rangle \subset W$ such that $\omega(w_1,w_2)\ne 0$ then $F = 0$.
\end{lemma}
\begin{proof}
Assume that $F \ne 0$. Let $q$ be the maximal integer such that $\CH^q(F)\ne 0$,
take a point $x\in\supp\CH^q(F)$ and choose $\langle w_1,w_2\rangle \subset W$
such that $x\in X_\phi$ (this is equivalent to the orthogonality of $w_1$ and $w_2$
with the 2-dimensional subspace of $W$ corresponding to $x\in X = \SGr(2,W)$).
Since the functor $i_\phi^*$ is left-exact it easily follows that $\CH^q(i_\phi^* F) \ne 0$, so $i_\phi^*F \ne 0$.
\end{proof}

Thus we have proved that the desired collection is indeed full.

\begin{remark}
The same argument allows to show that the Lefschetz exceptional collection
on a smooth hyperplane section of $\Gr(2,2m+1)$ obtained by proposition~\ref{lhp}
from the Lefschetz collection~\eqref{lcgr2mp1} is full.
\end{remark}

\section{Spinor bundles}\label{s_sb}

Our further goal is to construct an exceptional collection on the orthogonal isotropic
Grassmannian $\OGr(2,2m+1)$. However, in this case, as in the case of quadrics,
the restrictions of the tautological vector bundles from $\Gr(2,2m+1)$ don't give
a full exceptional collection and we need to consider some analogs of spinor bundles.
In this section we construct spinor bundles on isotropic Grassmannians
of a quadratic form, generalizing the definition of spinor bundles on quadrics
and investigate their properties.


We start with a reminder on Clifford algebras. Let $R$ be a commutative algebra over a field of zero characteristic,
$E$ a free $R$-module of rank $n$, and $q\in S^2 E^*$, a quadratic form.
The Clifford algebra of $q$ is defined as the quotient of the tensor algebra of $E$
by the following two-sided ideal (see \cite{Bou})
$$
\CB_q = T^\bullet(E)/\langle e_1\otimes e_2 + e_2\otimes e_1 - 2q(e_1,e_2){\mathbf{1}}\rangle.
$$
The Clifford algebra is naturally $\ZZ/2\ZZ$-graded, $\CB_q = \CB_q^0 \oplus \CB_q^1$,
the grading is induced by the $\ZZ$-grading of the tensor algebra. In what follows
we will be mostly interested in~$\CB_q^0$, the even part of the Clifford algebra.
Sometimes, instead of the quadratic form we will use the underlying vector space
as an index and sometimes (when the quadratic form and the space are clear)
the index will be omitted.
Note that as an $R$-module the Clifford algebra takes form
$$
\CB_q = \Lambda^\bullet E = R \oplus E \oplus \Lambda^2 E \oplus \dots,\qquad
\CB_q^0 = \Lambda^+ E = R \oplus \Lambda^2 E \oplus \Lambda^4 E \dots
$$

Assume that $n = 2m$ is even and the determinant of the quadratic form is invertible.
Then it is well known (\cite{Bou}) that if $q$ is neutral (i.e. the space $E$ of $q$
can be decomposed into a sum of two isotropic $R$-submodules) then $\CB^0$ is isomorphic
to a product of two matrix algebras.
Indeed, let $E = E_1 \oplus E_2$ be a decomposition of $E$ into a direct
sum of two isotropic submodules. Then the form $q$ gives an isomorphism
$E_2 \cong E_1^*$. Let $\BS_+$ and $\BS_-$ denote the even and the odd parts
of the exterior algebra of $E_1$ respectively:
$$
\BS_+ = \Lambda^+ E_1,\qquad
\BS_- = \Lambda^- E_1.
$$
The algebra $\CB$ acts on $\BS_+ \oplus \BS_-$ (the action of elements of $E_1$ is given
by the wedge-product, while the action of $E_2 \cong E_1^*$ is given by the convolution)
and the summands are invariant with respect to the action of $\CB^0$. Thus we get
a morphism of algebras $\CB^0 \to \End(\BS_+) \times \End(\BS_-)$ which is actually
an isomorphism. The simple $\CB^0$-modules $\BS_+$ and $\BS_-$ are called the {\sf half-spinor} modules.

The spinor module $\BS = \BS_+ \oplus \BS_-$ over $\CB$ is self-dual (nondegenerate $\CB$-invariant pairing
on $\BS \cong \Lambda^\bullet(E_1)$ is given by the determinant of the wedge-product). The duality isomorphism
$\BS_+ \oplus \BS_- \to \BS_+^* \oplus \BS_-^*$ gives a duality for half-spinor modules
(interchanging them if $m = n/2$ is odd).

Similarly, assume that $n = 2m+1$ is odd and the quadratic form is nondegenerate.
Let $E = E_1 \oplus E_2 \oplus Re$ be a decomposition of $E$ into a direct
sum of two isotropic submodules and orthogonal to them one-dimensional module.
The multiplication by $e$ gives a map $E_1\oplus E_2 \to \CB_q^0$ giving rise
to an isomorphism $\CB_{q'} \to \CB_q^0$, where $q'$ is a quadratic form on $E_1\oplus E_2$
given by the formula $q'(x) = - q(x)q(e)$. In particular, the spinor $\CB_{q'}$-module
$\BS = \Lambda^\bullet E_1$ acquires a structure of a $\CB_q^0$-module and it is known
that the corresponding homomorphism $\CB_q^0 \to \End(\BS)$ is an isomorpism.

In what follows we are going to consider the odd-dimensional and the even-dimensional cases
simultaneously (as far as it is possible). The discussed above differences between these cases
suggest to use the following convention. Let $\epsilon$ be an index taking values $+$, $-$, or empty,
with the following meaning. If $\epsilon=+$ (resp.\ $\epsilon = -$) this means that we are considering
the even-dimensional case and the object with this index corresponds to the half-spinor module $\BS_+$
(resp.\ $\BS_-$). On the other hand, if $\epsilon$ is empty this just means that we are considering
the odd-dimensional case. Having fixed this convention we can write down the above results as follows
$$
\CB^0 \cong \prod_{\epsilon} \End(\BS_\epsilon),
\qquad
\BS_\epsilon^*\cong \BS_{\pm\epsilon}.
$$

Below we will need also some results concerning the relation of a Clifford algebra (and its spinor modules)
for a quadratic form and some its subforms. Explicitly, assume that $U \subset E$ is an isotropic subspace, $\dim_R U = k$.
Let $U^\perp$ denote the orthogonal complement to $U$ in $E$, so that $U \subset U^\perp$.
Note that the initial nondegenerate quadratic form on $E$ induces a quadratic form $q'$ on $U^\perp \subset E$
with kernel $U$, and a nondegenerate quadratic form $q''$ on $U^\perp/U$. The embedding $U^\perp \subset E$
and the projection $U^\perp \to U^\perp/U$ are compatible with quadratic forms, hence induce morphisms
of Clifford algebras
$$
\xymatrix{
\CB_E^0 && \CB_{U^\perp/U}^0 \\
& \CB_{U^\perp}^0 \ar[ul] \ar[ur]
}
$$
Assume that we are given a decomposition $E = E_1 \oplus E_2$ if $n$ is even,
and a decomposition $E = E_1 \oplus E_2 \oplus Re$ if $n$ is odd,
where $E_1$, $E_2$ are isotropic $R$-submodules in $E$ and
$e \in E$ is orthogonal both to $E_1$ and $E_2$.
Assume also that $U \subset E_1$. Then we obtain also a decomposition
for $U^\perp/U$: $U^\perp/U = E_1/U \oplus E_2 \cap U^\perp$ if $n$ is even,
and $U^\perp/U = E_1/U \oplus E_2\cap U^\perp \oplus Re$ if $n$ is odd.
Let $\BS_{U,\epsilon}$ denote the (half)-spinor module of the algebra $\CB_{U^\perp/U}^0$.
We can consider both $\BS_\epsilon$ and $\BS_{U,\epsilon}$ as $\CB_{U^\perp}^0$-modules.

\begin{lemma}\label{spf1}
There is a canonical filtration $F^U_\bullet$ on the (half)-spinor $\CB_E^0$-module $\BS_\epsilon$ by $\CB_{U^\perp}^0$-submodules
$0 = F^U_{k+1}\BS_\epsilon \subset F^U_kS_\epsilon \subset \dots \subset F^U_1\BS_\epsilon \subset F^U_0\BS_\epsilon = \BS_\epsilon$
such that
$$
F^U_t\BS_\epsilon/F^U_{t+1}\BS_\epsilon \cong \BS_{U,(-1)^t\epsilon}\otimes\Lambda^tU.
$$
\end{lemma}
\begin{proof}
Actually, the Clifford algebra $\CB_{U^\perp/U}^0$ is isomorphic to the semisimple part of the algebra $\CB_{U^\perp}^0$
(the radical of $\CB_{U^\perp}^0$ equals $U\cdot\CB_{U^\perp}^1 \subset \CB_{U^\perp}^0$), and the desired filtration is
just the radical filtration of the $\CB_{U^\perp}^0$-module $\BS_\epsilon$, $F^U_t \BS_\epsilon = \Lambda^t U \cdot (\BS_{(-1)^t\epsilon})$.
The quotients of the radical filtration are modules over the semisimple part of the algebra,
hence they are isomorphic to direct sums of (half)-spinor modules. To see that they
have the form written up in the lemma
we note that the short exact sequence
$$
0 \to U \to E_1 \to E_1/U  \to 0
$$
induces a filtration on $\BS_\epsilon = \Lambda^\epsilon E_1$ with quotients
$\Lambda^{(-1)^t\epsilon}(E_1/U)\otimes\Lambda^t U = \BS_{U,(-1)^t\epsilon}\otimes\Lambda^tU$.
This filtration coincides with the discussed above radical filtration.
\end{proof}

\begin{lemma}\label{spf2}
Any choice of splitting $E_1 = U \oplus E_1'$ induces an isomorphism
of the Clifford algebra $\CB_{U^\perp}^0$ with $\prod_\epsilon \Hom(\BS_{U,\epsilon},\BS_\epsilon)$.
The isomorphism depends on a choice of splitting in such a way that the filtration on
$\prod_\epsilon \Hom(\BS_{U,\epsilon},\BS_\epsilon)$ induced by the filtration
$$
0 =
F_{\lfloor\frac{n-k}2\rfloor+1}\CB_{U^\perp}^0 \subset
F_{\lfloor\frac{n-k}2\rfloor}\CB_{U^\perp}^0 \subset
\dots \subset
F_1\CB_{U^\perp}^0 \subset
F_0\CB_{U^\perp}^0 =
\CB_{U^\perp}^0 =
\mathop{\oplus}\limits_{s \ge 0} \Lambda^{2s} U^\perp,
\qquad
F_t\CB_{U^\perp}^0 = \mathop{\oplus}\limits_{s \ge t} \Lambda^{2s} U^\perp
$$
on $\CB_{U^\perp}^0$ doesn't depend on a choice of splitting.
\end{lemma}
\begin{proof}
A choice of splitting
\begin{equation}\label{spl}
E_1 = U \oplus E_1'
\end{equation}
gives a decomposition
$U^\perp = U \oplus E'$, where
$E' = E'_1 \oplus E_2\cap U^\perp$ if $n$ is even, and
$E' = E'_1 \oplus E_2\cap U^\perp \oplus Re$ if $n$ is odd.
Note that the embeddings $U \subset U^\perp$ and $E' \subset U^\perp$
induce embeddings of the Clifford algebras $\Lambda^\bullet U = \CB_U \subset \CB_{U^\perp}$ and
$\CB_{E'} \subset \CB_{U^\perp}$. These subalgebras commute and we have an isomorphism
$\CB_{U^\perp} \cong \Lambda^\bullet U \otimes \CB_{E'}$, and
$\CB^0_{U^\perp} \cong (\Lambda^\bullet U \otimes \CB_{E'})^0$.

On the other hand, the splitting (\ref{spl}) gives a splitting of the filtration of lemma~\ref{spf1}
and furthermore an isomorphism $\BS_\epsilon \cong (\Lambda^\bullet U\otimes (\oplus_\epsilon \BS_{U,\epsilon}))^\epsilon$.
Thus we have
$$
\prod_\epsilon \Hom(\BS_{U,\epsilon},\BS_\epsilon) \cong
\prod_\epsilon \Hom(\BS_{U,\epsilon},(\Lambda^\bullet U\otimes (\oplus_\epsilon \BS_{U,\epsilon}))^\epsilon) \cong
(\Lambda^\bullet U\otimes \prod_\epsilon \Hom(\BS_{U,\epsilon},\BS_{U,\epsilon}))^0 =
(\Lambda^\bullet U \otimes \CB_{E'})^0
$$
and as a consequence we obtain an isomorphism
$\CB_{U^\perp}^0 \cong \prod_\epsilon \Hom(\BS_{U,\epsilon},\BS_\epsilon)$.
On the other hand, we have a canonical direct sum decomposition
$\CB_{U^\perp}^0 \cong \oplus \Lambda^{2s}U^\perp$, which via the constructed isomorphism
gives a direct sum decomposition of $\prod_\epsilon \Hom(\BS_{U,\epsilon},\BS_\epsilon)$.

However, the direct sum decomposition depends on a choice of splitting~(\ref{spl}).
So we are going to observe that the image of the filtration $F_t\CB_{U^\perp}^0$ on
$\prod_\epsilon \Hom(\BS_{U,\epsilon},\BS_\epsilon)$ doesn't depend on a choice of splitting.
Indeed, recall that the isomorphism
$\CB_{U^\perp}^0 \cong \prod_\epsilon \Hom(\BS_{U,\epsilon},\BS_\epsilon)$
which we use to transport the filtration is actually a composition of two isomorphisms
$$
\CB_{U^\perp}^0 \to (\Lambda^\bullet U \otimes \CB_{E'})^0 \to \prod_\epsilon \Hom(\BS_{U,\epsilon},\BS_\epsilon)
$$
It is clear that the first of these isomorphisms takes the direct sum decomposition of $\CB_{U^\perp}^0$
to the direct sum decomposition
$(\Lambda^\bullet U \otimes \CB_{E'})^0 =
(\Lambda^\bullet U \otimes \Lambda^\bullet E')^0 =
\oplus_s (\oplus_{k + l = 2s} \Lambda^{k} U \otimes \Lambda^{l} E')$,
so we are interested how the second isomorphism affects the latter direct sum decomposition.
Note that $\prod_\epsilon \Hom(\BS_{U,\epsilon},\BS_\epsilon)$ has a canonical structure of a right
$(\Lambda^\bullet U \otimes \CB_{U^\perp/U})^0$-module, and that the above isomorphism
is an isomorphism of right $(\Lambda^\bullet U \otimes \CB_{U^\perp/U})^0$-modules
(we have canonical isomorphism $E' \cong U^\perp/U$ compatible with quadratic forms,
showing that $(\Lambda^\bullet U \otimes \CB_{E'})^0$ is a free right
$(\Lambda^\bullet U \otimes \CB_{U^\perp/U})^0$-module of rank $1$).
It follows that a change of splitting~(\ref{spl}) results in an automorphism of
the $(\Lambda^\bullet U \otimes \CB_{U^\perp/U})^0$-module structure on
$\prod_\epsilon \Hom(\BS_{U,\epsilon},\BS_\epsilon)$, i.e. is given by the left action
of an invertible element of $(\Lambda^\bullet U \otimes \CB_{U^\perp/U})^0$.
A change of splitting can be written as $e' \mapsto e' + \phi(e')$, where $\phi:E'_1 \to U$
is an $R$-linear map. It is easy to check that the corresponding invertible element
of $(\Lambda^\bullet U \otimes \CB_{U^\perp/U})^0$ equals to
$$
(\Lambda^\bullet\phi)^0 \in
\Hom(\Lambda^\bullet E'_1,\Lambda^\bullet U)^0 =
(\Lambda^\bullet (E'_1)^* \otimes \Lambda^\bullet U)^0 \cong
(\Lambda^\bullet E'_2 \otimes \Lambda^\bullet U)^0 \subset
(\Lambda^\bullet U \otimes \CB_{E'})^0.
$$
Finally, we check that $(\Lambda^\bullet\phi)^0$
takes $\Lambda^{t_0} U \otimes \Lambda^{t_1} E'_1 \otimes \Lambda^{t_2} E'_2$
to $\oplus_{q_1,q_2\ge 0} \Lambda^{t_0+q_1+q_2} U \otimes \Lambda^{t_1-q_1} E'_1 \otimes \Lambda^{t_2+q_2} E'_2$,
and it is easy to see that the required filtration is preserved by this action.
\end{proof}

\bigskip

Now let $W$ be an orthogonal vector space (i.e. a vector space equipped with a nondegenerate quadratic form).
Put $n = \dim W$, $m=\lfloor n/2\rfloor$. Consider the isotropic Grassmannians $\OGr(m,W)$ of $m$-dimensional
(maximal) isotropic subspaces in $W$. It is well known that for odd $n$ the Grassmannian $\OGr(m,W)$ is connected,
its Picard group is $\ZZ$, its positive generator is very ample and the space of its global sections
is canonically isomorphic to the spinor module $\BS$ over the even part of the Clifford algebra $\CB_W^0$
(and the composition of embeddings $\OGr(m,W) \to \Gr(m,W) \subset \PP(\Lambda^mW)$ is given by the twice generator
of the Picard group).
Similarly, for even $n$ the Grassmannian $\OGr(m,W)$ has two connected components,
their Picard groups are $\ZZ$, positive generators are very ample and the spaces of their global sections
are canonically isomorphic to the half-spinor modules $\BS_\pm$ over the even part of the Clifford algebra $\CB_W^0$.
Let $\FF_m^\epsilon$ denote (the connected component of) the isotropic Grassmannian $\OGr(m,W)$
corresponding to the (half)-spinor module $\BS_\epsilon$, so that
$$
\begin{array}{ll}
\OGr(m,W) = \FF_m, & \text{if $n = 2m + 1$}\\
\OGr(m,W) = \FF_m^+ \sqcup \FF_m^-, & \text{if $n = 2m$}
\end{array}
\qquad\text{and}\qquad
H^0(\FF_m^\epsilon,\CO_{\FF_m^\epsilon}(1)) = \BS_\epsilon.
$$
where $\CO_{\FF_m^\epsilon}(1)$ is the positive generator of the Picard gorup of $\FF_m^\epsilon$.


Consider also the other isotropic Grassmannians in $W$, denote $\FF_k = \OGr(k,W)$.
Further, for any subset $I\subset\{1,2,\dots,m\}$ denote by $\FF_I^\epsilon$ the incidence
subvariety in (the connected component of) the product $\prod_{i\in I}\OGr(i,W)$ of the isotropic Grassmannians
and for $J \subset I$ denote by $\pi_J$ the projection $\FF^\epsilon_I \to \FF^\epsilon_J$.
Note that in this notation $\FF_1 = \OGr(1,W) = Q$ is the quadric in $\PP(W)$
and $\FF_2 = \OGr(2,W)$ is the isotropic Grassmannian of $2$-dimensional
subspaces which will be considered in the next section. Let also $\CO_{\FF^\epsilon_i}(1)$
denote the ample generator of the Picard group of $\FF^\epsilon_i$ and put
$$
\CO_{\FF^\epsilon_{i_1,i_2,\dots,i_k}}(d_1,d_2,\dots,d_k) :=
\pi_{i_1}^*\CO_{\FF^\epsilon_{i_1}}(d_1)\otimes
\pi_{i_2}^*\CO_{\FF^\epsilon_{i_1}}(d_2)\otimes\dots\otimes
\pi_{i_k}^*\CO_{\FF^\epsilon_{i_1}}(d_k).
$$
Let also $\CU_k$ denote the tautological rank $k$ subbundle in the trivial bundle $W\otimes\CO_{\FF_k}$,
the restriction of the tautological subbundle from $\Gr(k,W)$ to $\FF_k$.
Note that if $n$ is odd then $\FF_{1,2,\dots,m}$ is the flag variety of the group $\SO(W)$
and all $\FF_I$ are the partial flag varieties. Similarly, if $n$ is even then the flag variety of the group $\SO(W)$
is the fiber product $\FF^+_{1,2,\dots,m-2,m}\times_{\FF_{1,2,\dots,m-2}}\FF^-_{1,2,\dots,m-2,m}$
and all $\FF^\epsilon_I$ with $m-1\not\in I$ are the partial flag varieties.

\smallskip

Now take any $I \subset\{1,2,\dots,m-1\}$, consider the diagram
$$
\xymatrix{
& \FF^\epsilon_{I,m} \ar[dl]_{\pi_I} \ar[dr]^{\pi_m} \\
\FF_I && \FF^\epsilon_m
}
$$
and define the spinor bundles on $\FF_I$ by the formula
\begin{equation}\label{spinor}
\CS_{I\epsilon} = \pi_{I*}\pi_m^*\CO_{\FF^\epsilon_m}(1).
\end{equation}
Note that taking $I = \{1\}$ we obtain the usual spinor bundles on the quadric $Q = \FF_1$ (see \cite{Ka}),
and taking $I = \emptyset$ we obtain the (half)-spinor modules $\BS_\epsilon$
(considered as vector bundles on $\FF_\emptyset = \Spec k$).
Note also that
\begin{equation}\label{spinorI}
\CS_{I\epsilon} \cong \pi_k^*\CS_{k\epsilon},
\qquad\text{where $k = \max\{i\in I\}$}.
\end{equation}
Indeed, this easily follows from the base change since $\FF^\epsilon_{I,m} = \FF_I\times_{\FF_k} \FF^\epsilon_{k,m}$.

It is important that the spinor bundles can be considered as a generalization of spinor modules in a relative situation.
Indeed, consider the vector bundle $\CU_k^\perp/\CU_k$ over $\FF_k$. It carries a natural quadratic form induced by the
quadratic form on $W$. Let $\CB_k^0$ denote the corresponding sheaf of even parts of Clifford algebras  on $\FF_k$.
The fiber of $\CB_k^0$ over a point of $\FF_k$ corresponding to an isotropic subspace $U \subset W$ is the Clifford algebra
$\CB_{U^\perp/U}^0$ and the fiber of $\CS_{k\epsilon}$ is the (half)-spinor module $\BS_{U,\epsilon}$. So we can use
local properties of spinor modules proved in lemma~\ref{spf1} and \ref{spf2} to deduce the global properties
of spinor bundles.

\begin{proposition}\label{filtcs}
The pullback of the (half)-spinor module $\BS_\epsilon$ to $\FF_k$ admits a length $k+1$ filtration
$0 = F_{k+1}\BS_\epsilon \subset F_k \BS_\epsilon \subset \dots \subset F_1\BS_\epsilon \subset F_0\BS_\epsilon = \BS_\epsilon\otimes\CO_{\FF_k}$
such that
$$
F_t\BS_\epsilon/F_{t+1}\BS_\epsilon \cong \CS_{k,(-1)^t\epsilon}\otimes\Lambda^t\CU_k.
$$
In particular, there is a short complex on $\FF_2 = \OGr(2,W)$ of the form
\begin{equation}\label{csu}
0 \to \CS_{2,\epsilon}\otimes\CO_X(-1) \to \BS_\epsilon\otimes\CO_X \to \CS_{2,\epsilon} \to 0
\end{equation}
with the only nontrivial cohomology in the middle isomorphic to $\CS_{2,-\epsilon}\otimes\CU_2$.
\end{proposition}
\begin{proof}
This is a global version of the local filtration of lemma~\ref{spf1}.
\end{proof}

Similarly, one can prove a relative version of this proposition.

\begin{proposition}\label{filtsk}
The pullback $\pi_k^*\CS_{k\epsilon}$ of the spinor bundle to $\FF_{k,l}$ admits a length $(l-k+1)$ filtration
$0 = F_{l-k+1}\pi_k^*\CS_{k\epsilon} \subset F_{l-k}\pi_k^*\CS_{k\epsilon} \subset \dots \subset F_1\pi_k^*\CS_{k\epsilon} \subset F_0\pi_k^*\CS_{k\epsilon} = \pi_k^*\CS_{k\epsilon}$
such that
$$
F_t\pi_k^*\CS_{k\epsilon}/F_{t+1}\pi_k^*\CS_{k\epsilon} \cong \CS_{l,(-1)^t\epsilon}\otimes\Lambda^t(\CU_l/\CU_k).
$$
In particular, we have the following short exact sequence on $\FF_{k-1,k}$
\begin{equation}\label{csk}
0 \to \pi_k^*\CS_{k,-\epsilon}\otimes\CO_{\FF_{k-1,k}}(1,-1) \to \pi_{k-1}^*\CS_{k-1,\epsilon} \to \pi_k^*\CS_{k\epsilon} \to 0.
\end{equation}
\end{proposition}

\begin{corollary}\label{rd}
We have
$$
\rank(\CS_{k\epsilon}) =
\begin{cases}
2^{m-k}, & \text{if $n=2m+1$}\\
2^{m-k-1}, & \text{if $n=2m$}
\end{cases}
\qquad
\det\CS_{k\epsilon} \cong
\begin{cases}
\CO_{\FF_k}(2^{m-k-1}), & \text{if $n=2m+1$}\\
\CO_{\FF_k}(2^{m-k-2}), & \text{if $n=2m$}
\end{cases}
$$
\end{corollary}
\begin{proof}
For $k=0$ the claim is evident. For $k>0$ we proceed by induction using~(\ref{csk}).
\end{proof}

\begin{proposition}\label{selfd}
The spinor bundles are (mutually) selfdual up to a twist: $\CS^*_{k\epsilon} \cong \CS_{k,(-1)^m\epsilon}(-1)$.
\end{proposition}
\begin{proof}
We have seen that the (half)-spinor modules over the Clifford algebras are (mutually) self-dual
and the duality isomorphisms (compatible with the module structures) is unique, hence the local
isomorphisms glue into an isomorphism $\CS^*_{k\epsilon} \cong \CS_{k,(-1)^m\epsilon}\otimes L$
for some line bundle $L$ on $\FF_k$. But comparing determinants and using corollary~\ref{rd},
we easily deduce that $L \cong \CO_{\FF_k}(-1)$.
\end{proof}

\begin{proposition}\label{csf2}
There is a filtration on the bundle $\oplus_\epsilon \BS_\epsilon\otimes \CS_{k\epsilon}^*$
$$
0 =
F_{\lfloor\frac{n-k}2\rfloor+1}(\oplus_\epsilon \BS_\epsilon\otimes \CS_{k\epsilon}^*) \subset
F_{\lfloor\frac{n-k}2\rfloor}(\oplus_\epsilon \BS_\epsilon\otimes \CS_{k\epsilon}^*) \subset
\dots \subset
F_1(\oplus_\epsilon \BS_\epsilon\otimes \CS_{k\epsilon}^*) \subset
F_0(\oplus_\epsilon \BS_\epsilon\otimes \CS_{k\epsilon}^*) =
\oplus_\epsilon \BS_\epsilon\otimes \CS_{k\epsilon}^*,
$$
such that
$F_s(\oplus_\epsilon \BS_\epsilon\otimes \CS_{k\epsilon}^*)/F_{s+1}(\oplus_\epsilon \BS_\epsilon\otimes \CS_{k\epsilon}^*) \cong \Lambda^{2s} \CU_k^\perp$.
\end{proposition}
\begin{proof}
This is a global version of the local filtration of lemma~\ref{spf2}.
\end{proof}

\begin{proposition}\label{skexc}
The spinor bundles on $\FF_k$ are exceptional and $H^\bullet(\FF_k,\CS_{k\epsilon}(-p)) = 0$ for $1\le p\le n-2k$.
Moreover, if $2\le k\le m-1$ then the collection
$$
\begin{array}{ll}
(\CS_k,\CS_k(1),\dots,\CS_k(n - 2k)), & \text{if $n = 2m + 1$}\\
(\CS_{k-},\CS_{k+},\CS_{k-}(1),\CS_{k+}(1),\dots,\CS_{k-}(n - 2k),\CS_{k+}(n - 2k)), & \text{if $n = 2m$}
\end{array}
$$
is exceptional.
\end{proposition}
\begin{proof}
In fact, all this can be checked via the Borel--Bott--Weil theorem.
However, this approach requires the combinatorics of the Weil group action
on the weight lattice, and is not very illuminating. So, we give here another proof,
using exceptionality of spinor bundles on quadrics (proved by Kapranov \cite{Ka})
and exact sequences~(\ref{csk}).

Let us compute $\RHom(\CS_{k\epsilon}(p),\CS_{k\epsilon'})$.
Consider the diagram
$$
\xymatrix{
& \FF_{k-1,k} \ar[dl]_{\pi_{k-1}} \ar[dr]^{\pi_k} \\
\FF_{k-1} && \FF_k
}
$$
Note that the map $\pi_k$ is a projectivization of a vector bundle, hence the functor $\pi_k^*$
is fully faithful, i.e. we can (and will) compute
\begin{multline*}
\RHom(\pi_k^*\CS_{k\epsilon}(p),\pi_k^*\CS_{k\epsilon'}) =
\RHom(\CS_{k-1,k\epsilon}(0,p),\CS_{k-1,k\epsilon'}) = \\
\RGamma(\FF_{k-1,k},\CS_{k-1,k\epsilon}^*\otimes\CS_{k-1,k\epsilon'}(0,-p)) =
\RGamma(\FF_{k-1},\pi_{k-1*}(\CS_{k-1,k\epsilon}^*\otimes\CS_{k-1,k\epsilon'}(0,-p)))
\end{multline*}
instead. Now look at the map $\pi_{k-1}$. It is a fibration in quadrics
(of dimension $n-2(k-1)-2 = n-2k$) and it is clear that restriction
of the spinor bundle $\CS_{k-1,k\epsilon}$ to any fiber of $\pi_{k-1}$
is isomorphic to the usual spinor bundle on this fiber. Therefore
on any fiber of $\pi_{k-1}$ the bundle $\CS_{k-1,k\epsilon}^*\otimes\CS_{k-1,k\epsilon'}$
has no cohomology if $\epsilon\ne\epsilon'$ and has one-dimensional cohomology in degree 0
if $\epsilon = \epsilon'$. Since $\pi_{k-1}$ is flat, it follows that
$$
\pi_{k-1*}(\CS_{k-1,k\epsilon}^*\otimes\CS_{k-1,k\epsilon'}) =
\begin{cases}
0, & \text{if $\epsilon \ne \epsilon'$}\\
\text{a line bundle}, & \text{if $\epsilon = \epsilon'$}
\end{cases}
$$
On the other hand, if $\epsilon = \epsilon'$, the bundle
$\CS_{k-1,k\epsilon}^*\otimes\CS_{k-1,k\epsilon'}$ contains $\CO_{\FF_{k-1,k}}$ as a direct summand,
hence its push-forward to $\FF_{k-1}$ contains $\pi_{k-1*}(\CO_{\FF_{k-1,k}}) \cong \CO_{\FF_{k-1}}$
as a direct summand, hence the line bundle discussed above is trivial. Summarizing, we have
$$
\pi_{k-1*}(\CS_{k-1,k\epsilon}^*\otimes\CS_{k-1,k\epsilon'}) =
\begin{cases}
0, & \text{if $\epsilon \ne \epsilon'$}\\
\CO_{\FF_{k-1}}, & \text{if $\epsilon = \epsilon'$}
\end{cases}
$$
Similarly, for $0\le p\le n-2k-1$ the bundle $\CS_{k-1,k\epsilon}^*(0,-p)$ has no cohomology on any fiber of $\pi_{k-1}$, hence
\begin{equation}\label{piscsp}
\pi_{k-1*}(\CS_{k-1,k\epsilon}^*(0,-p)) = 0
\qquad\text{for $0\le p\le n-2k-1$ and any $\epsilon$.}
\end{equation}
Now consider the sequence~(\ref{csk}) tensored by $\CS_{k-1,k\epsilon}^*(0,-p)$
$$
0 \to
\CS_{k-1,k\epsilon}^*\otimes\CS_{k-1,k,-\epsilon'}(1,-p-1) \to
\CS_{k-1,k\epsilon}^*(0,-p)\otimes\pi_{k-1}^*\CS_{k-1,\epsilon'} \to
\CS_{k-1,k\epsilon}^*\otimes\CS_{k-1,k\epsilon'}(0,-p) \to
0.
$$
Taking the push-forward to $\FF_{k-1}$ and using~(\ref{piscsp}) and the projection formula, we deduce that
$$
\pi_{k-1*}(\CS_{k-1,k\epsilon}^*\otimes\CS_{k-1,k,-\epsilon'}(0,-p-1)) \cong
\pi_{k-1*}(\CS_{k-1,k\epsilon}^*\otimes\CS_{k-1,k\epsilon'}(0,-p))\otimes\CO_{\FF_{k-1}}(-1)[-1],
$$
for $0\le p\le n-2k-1$. It follows by induction that
$$
\pi_{k-1*}(\CS_{k-1,k\epsilon}^*\otimes\CS_{k-1,k\epsilon'}(0,-p)) =
\begin{cases}
0, & \text{if $\epsilon \ne (-1)^p\epsilon'$, $0\le p\le n-2k$}\\
\CO_{\FF_{k-1}}(-p)[-p], & \text{if $\epsilon = (-1)^p\epsilon'$, $0\le p\le n-2k$}
\end{cases}
$$
It remains to note that $\FF_{k-1}$ is a Fano variety with canonical bundle isomorphic to
$\CO_{\FF_{k-1}}(k-n)$ therefore the bundles $\CO_{\FF_{k-1}}(-p)$ have no cohomology on $\FF_{k-1}$
for $1\le p\le n-2k$, hence
$$
\RHom(\pi_k^*\CS_{k\epsilon}(p),\pi_k^*\CS_{k\epsilon'}) =
\begin{cases}
0, & \text{if $\epsilon \ne \epsilon'$, $0\le p\le n-2k$, or $\epsilon = \epsilon'$ and $1\le p\le n-2k$}\\
\kk, & \text{if $\epsilon = \epsilon'$ and $p=0$}
\end{cases}
$$
and this is precisely what we need.

It remains to note that $\CS_{k\epsilon}(-1) \cong \CS_{k,-\epsilon}^*$ by proposition~\ref{selfd},
so~(\ref{piscsp}) implies the desired vanishing of cohomology of twists of $\CS_{k\epsilon}$.
\end{proof}

\section{Odd-dimensional orthogonal case}\label{s_ogr2mp1}

Consider the isotropic Grassmannian $X = \OGr(2,W)$ of two-dimensional isotropic subspaces in
an orthogonal vector space $W$ of dimension $n = 2m + 1$.
Let $\CU$ denote the restriction of the tautological rank $2$ subbundle from $\Gr(2,W)$ to $X = \OGr(2,W)$.
Let also $\CS = \CS_2$ denote the spinor bundle on $X$ constructed in the previous section.

Consider the Lefschetz collection~\eqref{lcogr2mp1} on $\OGr(2,W)$.
This is a Lefschetz collection with the first block
\begin{equation}\label{ugco}
(\CO_X,\CU^*,S^2\CU^*,\dots,S^{m-2}\CU^*,\CS)
\end{equation}
and with the support partition
\begin{equation}\label{ugpo}
\lambda^{\mathsf O} = (\underbrace{m,m,\dots,m}_{2m-2}).
\end{equation}
In the other words, the collection consists of vector bundles $S^l\CU^*(k)$ with integers $(k,l)$
from the set
\begin{equation}\label{Sno}
\Upso_{2m+1} = \{(k,l)\in\ZZ^2\ |\ \text{$0\le k\le 2m-3$, $0\le l\le m-2$}\}.
\end{equation}
and of an additional collection
\begin{equation}\label{spino}
(\CS,\CS(1),\dots,\CS(2m-3))
\end{equation}
of twists of the spinor bundle.
On the following picture we draw the regions $\Ups_{2m+1}$ and $\Upso_{2m+1}$
in the $(k,l)$-plane.
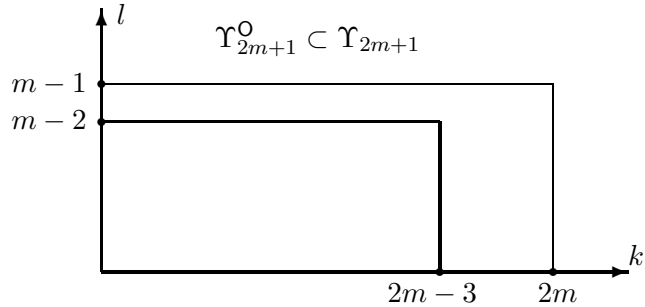
\begin{figure}[h]
\unitlength=1mm
\begin{picture}(170,45)
\thicklines
\put(75,35){$\Upso_{2m+1} \subset \Ups_{2m+1}$}
\put(60,5){\vector(0,1){35}}
\put(60,5){\vector(1,0){70}}
\put(60,25){\line(1,0){45}}
\put(105,5){\line(0,1){20}}

\thinlines

\put(60,30){\line(1,0){60}}
\put(120,5){\line(0,1){25}}

\put(60,30){\circle*{1}}
\put(48,29){$m-1$}
\put(60,25){\circle*{1}}
\put(48,24){$m-2$}
\put(105,5){\circle*{1}}
\put(98,1){$2m-3$}
\put(120,5){\circle*{1}}
\put(118,1){$2m$}
\put(130,6){$k$}
\put(62,38){$l$}

%

\end{picture}
\end{figure}


\begin{theorem}
Let $\Upso_{2m+1}$ be the set defined in~\eqref{Sno}.
Then the Lefschetz collection $\{S^l\CU^*(k)\ |\ (k,l)\in\Upso_{2m+1}\} \cup \{\CS(k)\ |\ 0\le k\le 2m-3\}$
is a full exceptional collection in $\D^b(\OGr(2,2m+1))$.
\end{theorem}

First of all we must check that the collection is exceptional.

\begin{lemma}
The Lefschetz collection $\{S^l\CU^*(k)\ |\ (k,l)\in\Upso_{2m+1}\} \cup \{\CS(k)\ |\ 0\le k\le 2m-3\}$ is exceptional.
\end{lemma}
\begin{proof}
First, let us check that the collection $\{S^l\CU^*(k)\ |\ (k,l)\in\Upso_{2m+1}\}$ is exceptional.
The arguments of lemma~\ref{critl} and lemma~\ref{bbwgr2n} show that it suffices to check that
$H^\bullet(X,\Sigma^{p-k,-l-k}\CU^*) = 0$ for $(k,l)\in\Upso_{2m+1}$ and $0\le p\le m-2$
(with an additional restriction $p<l$ for $k=0$), and that $H^\bullet(X,\CO_X) = \kk$.

For this we note that $X \subset \Gr(2,W)$ is the zero locus of a regular section of the bundle
$S^2\CU^*$ (given by the quadratic form on $W$). Therefore, on $\Gr(2,W)$ we have the following Koszul resolution
$$
0 \to \CO(-3) \to S^2\CU(-1) \to S^2\CU \to \CO \to i_*\CO_X \to 0.
$$
Tensoring it by $\Sigma^{p-k,-l-k}\CU^*$ we obtain a resolution of
$\Sigma^{p-k,-l-k}\CU^* \otimes i_*\CO_X \cong i_*(\Sigma^{p-k,-l-k}\CU^*)$
on $\Gr(2,W)$. It is easy to see that
$$
\begin{array}{lllllll}
\Sigma^{p-k,-l-k}\CU^*\otimes S^2\CU     &\cong& \Sigma^{p-k,-l-k-2}\CU^*   &\oplus& \Sigma^{p-k-1,-l-k-1}\CU^* &\oplus& [\Sigma^{p-k-2,-l-k}\CU^*],\\
\Sigma^{p-k,-l-k}\CU^*\otimes S^2\CU(-1) &\cong& \Sigma^{p-k-1,-l-k-3}\CU^* &\oplus& \Sigma^{p-k-2,-l-k-2}\CU^* &\oplus& [\Sigma^{p-k-3,-l-k-1}\CU^*],\\
\Sigma^{p-k,-l-k}\CU^*\otimes \CO(-3)    &\cong& \Sigma^{p-k-3,-l-k-3}\CU^*
\end{array}
$$
(the terms in the brackets should be omitted if $p+l\le 1$).
All summands in the RHS can be rewritten in the form
$\Sigma^{p'-k',-l'-k'}\CU^*$ with $(k',l') \in \Ups_{2m+1}$ and $0\le p\le m-1$
(with an additional restriction $p'<l'$ for~$k'=0$). As it was checked in lemma~\ref{bbwgr2n},
the cohomology of all these bundles on $\Gr(2,W)$ vanishes,
with the only exception $H^\bullet(\Gr(2,W),\CO_{\Gr(2,W)}) = \kk$,
therefore, the considered collection of vector bundles
on~$X$ is indeed exceptional.

Further, we note that the collection $(\CS,\CS(1),\dots,\CS(2m-3))$ is exceptional by proposition~\ref{skexc}.
It remains to check semiorthogonalities between $S^l\CU^*(k)$ and $\CS(t)$.

Now let us check that $\RHom(S^l\CU^*(k),\CS) = H^\bullet(X,\CS\otimes S^l\CU(-k)) = 0$ for all $(k,l)\in\Upso_{2m+1}$ with $k\ge 1$,
i.e. for all $1\le k\le 2m-3$, $0\le l\le m-2$. Actually, we will check that
\begin{equation}\label{hsu}
H^\bullet(X,\CS\otimes S^l\CU(-k)) = 0 \qquad \text{for all $0\le l\le m-2$, $1-l\le k\le 2m-2$.}
\end{equation}
The proof is inductive in $l$. The case $l=0$ is already proved in proposition~\ref{skexc}. Now consider the complex~(\ref{csu})
tensored by $S^{l-1}\CU(-k)$:
$$
0 \to
\CS\otimes S^{l-1}\CU(-k-1) \to
\BS\otimes S^{l-1}\CU(-k) \to
\CS\otimes S^{l-1}\CU(-k) \to 0
$$
By proposition~\ref{filtcs} its cohomology is isomorphic to $\CS\otimes \CU\otimes S^{l-1}\CU(-k)$.
Consider the hypercohomology spectral sequence of this complex. By the induction hypothesis and the proved above
exceptionality of the collection $\{S^l\CU^*(k)\ |\ (k,l)\in\Upso_{2m+1}\}$, the cohomology of all terms of this
complex vanish if $2-l\le k\le 2m-3$. Therefore, the cohomology of the bundle
$\CS\otimes \CU\otimes S^{l-1}\CU(-k)$ vanishes as well. But
$\CS\otimes S^{l}\CU(-k)$ is a direct summand of $\CS\otimes \CU\otimes S^{l-1}\CU(-k)$,
hence its cohomology vanishes as well. It remains to consider only the bundle $\CS\otimes S^{l}\CU(-k)$
with $k = 2m-2$ and $k=1-l$ which are not covered by the above arguments. Let us start with $k = 1-l$.
Note that $S^{l-1}\CU(l-1)\cong S^{l-1}\CU^*$, so the above complex takes form
$$
0 \to
\CS\otimes S^{l-1}\CU^*(-1) \to
\BS\otimes S^{l-1}\CU^* \to
\CS\otimes S^{l-1}\CU^* \to 0
$$
and its cohomology is isomorphic to
$\CS\otimes\CU\otimes S^{l-1}\CU^*\cong
\CS\otimes S^{l-2}\CU^* \oplus \CS\otimes S^{l}\CU^*(-1)$.
The cohomology of the first term of the complex vanishes by the induction hypothesis,
while for the middle and the last term of the complex the Borel--Bott--Weil theorem gives
$$
H^\bullet(X,\BS\otimes S^{l-1}\CU^*) \cong V(\omega_m) \otimes V((l-1)\omega_1),\qquad
H^\bullet(X,\CS\otimes S^{l-1}\CU^*) \cong V(\omega_m + (l-1)\omega_1),
$$
(the application of the Borel--Bott--Weil theorem in this case is absolutely straightforward
since the corresponding weights are dominant), and the map between them induced by
the differential of the complex is the canonical projection.
But according to the representation theory of the Lie group $\SO(2m+1)$ we have
$V(\omega_m) \otimes V((l-1)\omega_1) \cong V(\omega_m + (l-1)\omega_1) \oplus V(\omega_m + (l-2)\omega_1)$,
hence
$$
H^\bullet(X,\CS\otimes S^{l-2}\CU^* \oplus \CS\otimes S^{l}\CU^*(-1)) \cong V(\omega_m + (l-2)\omega_1).
$$
It remains to note that $H^\bullet(X,\CS\otimes S^{l-2}\CU^*) \cong V(\omega_m + (l-2)\omega_1)$
again by Borel--Bott-Weil, hence the cohomology of $\CS\otimes S^{l}\CU^*(-1) \cong \CS\otimes S^{l}\CU(l-1)$
vanishes, but this is precisely what we need.

Finally, the case $k=2m-2$ follows by Serre duality.
The canonical bundle of $X$ is isomorphic to $\CO_X(2-2m)$, hence
$$
H^{\dim X - \bullet}(X,\CS\otimes S^{l}\CU(2-2m)) \cong
H^\bullet(X,\CS^*\otimes S^l\CU^*))^*,
$$
but by proposition~\ref{selfd} $\CS^* \cong \CS(-1)$, while $S^l\CU^* \cong S^l\CU(l)$, hence
$$
\CS^*\otimes S^l\CU^* \cong \CS\otimes S^{l}\CU(l-1).
$$

Finally we note that by Serre duality we have
$\RHom(\CS(k),S^l\CU^*) \cong \RHom(S^l\CU^*,\CS(k+2-2m)) = H^\bullet(X,\CS\otimes S^l\CU(k+2-2m))$
which is zero for $0\le k\le 2m-3$ by (\ref{hsu}).
\end{proof}

It remains to prove the fullness of the collection.
We use the same method as in sections~\ref{s_gr} and~\ref{s_sgr}
with the induction step that changes $\dim W$ by $2$.

Consider the set
\begin{equation}\label{TSno}
\TUpso_{2m-1} = \{(k,l)\in\ZZ^2\ |\ \text{$0\le k\le 2m-3$, $0\le l\le m-1$}\}.
\end{equation}
On the following picture the region $\Upso_{2m+1}$ is drawn together with the region $\TUpso_{2m-1}$.

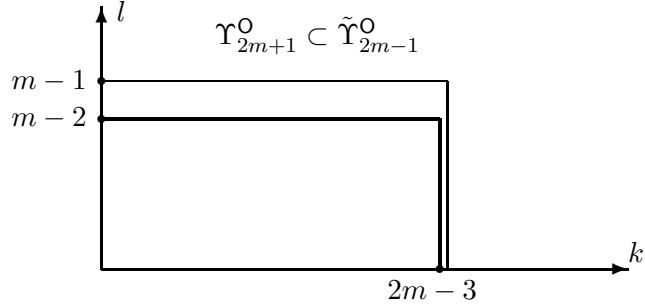
\begin{figure}[h]
\unitlength=1mm
\begin{picture}(170,45)
\thicklines
\put(75,35){$\Upso_{2m+1} \subset \TUpso_{2m-1}$}
\put(60,5){\vector(0,1){35}}
\put(60,5){\vector(1,0){70}}
\put(60,25){\line(1,0){45}}
\put(105,5){\line(0,1){20}}

\thinlines

\put(60,30){\line(1,0){46}}
\put(106,5){\line(0,1){25}}

\put(60,30){\circle*{1}}
\put(48,29){$m-1$}
\put(60,25){\circle*{1}}
\put(48,24){$m-2$}
\put(105,5){\circle*{1}}
\put(98,1){$2m-3$}
\put(130,6){$k$}
\put(62,38){$l$}


\end{picture}
\end{figure}


\begin{lemma}\label{tsigmao}
For any $(k,l)\in\TUpso_{2m-1}$ the vector bundle $S^l\CU^*(k)$ lies in the triangulated category
generated by the Lefschetz collection $\{S^l\CU^*(k)\ |\ (k,l)\in\Upso_{2m+1}\} \cup \{\CS(k)\ |\ 0\le k\le 2m-3\}$ in $\D^b(\OGr(2,W))$.
%
\end{lemma}
\begin{proof}
Let $\D$ denote the triangulated subcategory of $\D^b(\OGr(2,W))$ generated by the Lefschetz collection
$\{S^l\CU^*(k)\ |\ (k,l)\in\Upso_{2m+1}\} \cup \{\CS(k)\ |\ 0\le k\le 2m-3\}$. Note that
$$
\TUpso_{2m-1} \setminus \Upso_{2m+1} = \{(0,m-1),(1,m-1),\dots,(2m-3,m-1)\}.
$$
So, we have to check that the vector bundles
$$
S^{m-1}\CU^*,S^{m-1}\CU^*(1),\dots,S^{m-1}\CU^*(2m-3),
$$
lie in $\D$.
Consider the filtration on the sheaf $\BS^*\otimes\CS$, dual to the filtration
described in proposition~\ref{csf2}. The quotients of this filtration take form
$\Lambda^{2s}(W/\CU)$. Note that if $2s < m-1$ then $\Lambda^{2s}(W/\CU)$ admits
a resolution~(\ref{sku}), which shows that we have
$\Lambda^{2s}(W/\CU) \in \langle \CO_X,\CU^*(-1),S^2\CU^*(-2),\dots,S^{2s}\CU^*(-2s)\rangle$
hence
\begin{equation}\label{inc1}
\Lambda^{2s}(W/\CU)(t) \in \D
\qquad\text{if $2s < m-1$ and $2s\le t\le 2m-3$.}
\end{equation}
Similarly, if $2s > m$ then $\Lambda^{2s}(W/\CU) \cong \Lambda^{2m-1-2s}\CU^\perp$ admits
a resolution~(\ref{skus}), which shows that we have
$\Lambda^{2s}(W/\CU) \in \langle \CO_X(1),\CU^*(1),S^2\CU^*(1),\dots,S^{2m-1-2s}\CU^*(1)\rangle$
hence
\begin{equation}\label{inc2}
\Lambda^{2s}(W/\CU)(t) \in \D
\qquad\text{if $2s > m$ and $-1\le t\le 2m-4$.}
\end{equation}
But there exists precisely one even integer $2s$ which is not covered by the previous two cases.
It is either $m-1$ if $m$ is odd, or $m$ if $m$ is even. We are forced to consider this two cases
separately.

First assume that $m$ is odd. Then the inclusions~(\ref{inc1}), (\ref{inc2}) together with
\begin{equation}\label{inc3}
\BS^*\otimes\CS(t) \in \D
\qquad\text{if $0\le t\le 2m-3$}
\end{equation}
imply that for $m-3 \le t\le 2m-4$ both $\BS^*\otimes\CS(t)$ and all quotients of its filtration
except for $\Lambda^{m-1}(W/\CU)(t)$ lie in $\D$. But since $\D$ is triangulated it follows also that
$\Lambda^{m-1}(W/\CU)(t) \in \D$ for $m-3 \le t\le 2m-4$ as well.
But looking at the resolution~(\ref{sku}) we see that this implies that $S^{m-1}\CU^*(t-m+1) \in \D$
for $m-2 \le t \le 2m-4$. Thus we deduce that for odd $m$ we have
\begin{equation}\label{cinc1}
S^{m-1}\CU^*(-1),S^{m-1}\CU^*,\dots,S^{m-1}\CU^*(m-3) \in \D.
\end{equation}

Now assume that $m$ is even. Then the inclusions~(\ref{inc1}), (\ref{inc2}) and (\ref{inc3})
imply that for $m-2 \le t\le 2m-4$ both $\BS^*\otimes\CS(t)$ and all quotients of its filtration
except for $\Lambda^{m}(W/\CU)(t)$ lie in $\D$. But since $\D$ is triangulated it follows also that
$\Lambda^{m}(W/\CU)(t) \cong \Lambda^{m-1}\CU^\perp(t+1) \in \D$ for $m-2 \le t\le 2m-4$ as well.
But looking at the resolution~(\ref{skus}) we see that this implies that $S^{m-1}\CU^*(t+1) \in \D$
for $m-2 \le t \le 2m-4$. Thus we deduce that for even $m$ we have
\begin{equation}\label{cinc2}
S^{m-1}\CU^*(m-1),S^{m-1}\CU^*(m),\dots,S^{m-1}\CU^*(2m-3) \in \D.
\end{equation}

Finally, we are going to use the following trick to show that the other twists of $S^{m-1}\CU^*$
also lie in $\D$. Consider the vector bundle $\CU^\perp/\CU$ on $X$. It comes with a nondegenerate
quadratic form induced by the form on $W$. Its rank is $2m+1 - 4 = 2m-3$.
It follows that for any $t$ we have an isomorphism
\begin{equation}\label{liso}
\Lambda^{m-1}(\CU^\perp/\CU)(t) \cong \Lambda^{m-2}(\CU^\perp/\CU)(t).
\end{equation}
Now note that $\CU^\perp/\CU$ is the cohomology of the complex $\CU \to W\otimes\CO_X \to \CU^*$.
Taking exterior powers we deduce that
$$
\Lambda^s(\CU^\perp/\CU) \cong
\left\{\vcenter{\scriptsize
\xymatrix@=10pt{
S^s\CU \ar[r] & W\otimes S^{s-1}\CU \ar[r] \ar[d] & \dots \ar[r] & \Lambda^{s-2}W\otimes S^2\CU \ar[r] \ar[d] & \Lambda^{s-1}W\otimes \CU \ar[r] \ar[d] & \Lambda^sW\otimes\CO_X \ar[d] \\
& \CU^*\otimes S^{s-1}\CU \ar[r] & \dots \ar[r] & \CU^*\otimes\Lambda^{s-3}W\otimes S^2\CU \ar[r] \ar[d] & \CU^*\otimes\Lambda^{s-2}W\otimes \CU \ar[r] \ar[d] & \CU^*\otimes\Lambda^{s-1}W \ar[d] \\
&& \ddots & \vdots \ar[d] & \vdots \ar[d] & \vdots \ar[d] \\
&&& S^{s-2}\CU^*\otimes S^2\CU \ar[r] & S^{s-2}\CU^*\otimes W\otimes \CU \ar[r] \ar[d] & S^{s-2}\CU^* \otimes \Lambda^2W \ar[d] \\
&&&& S^{s-1}\CU^*\otimes \CU \ar[r] & S^{s-1}\CU^*\otimes W \ar[d]\\
&&&&& S^s\CU^*
}}\right\}
$$
Taking into account isomorphisms
$$
S^a\CU^*\otimes S^b\CU \cong \bigoplus_{c=0}^{\min\{a,b\}} S^{a+b-2c}\CU^*(c-b)
$$
we see that the isomorphism~(\ref{liso}) with $t=m-2$ can be considered as an expression for $S^{m-1}\CU^*(m-2)$
in terms of $S^{m-1}\CU^*(-1),\dots,S^{m-1}\CU^*(m-3)$ and other objects of $\D$. Taking into account~(\ref{cinc1})
we deduce that $S^{m-1}\CU^*(m-2)\in\D$ for odd $m$. Repeating the same argument for the isomorphisms~(\ref{liso})
with $t = m-1,m,\dots,2m-3$ we deduce that $S^{m-1}\CU^*(t) \in \D$ for these $t$ for odd $m$.
Similarly, if $m$ is even we see the isomorphism~(\ref{liso}) with $t=2m-3$ can be considered as an expression
for $S^{m-1}\CU^*(m-2)$ in terms of $S^{m-1}\CU^*(m-1),\dots,S^{m-1}\CU^*(2m-3)$ and other objects of $\D$.
Taking into account~(\ref{cinc2}) we deduce that $S^{m-1}\CU^*(m-2)\in\D$ for even $m$. Repeating the same argument
for the isomorphisms~(\ref{liso}) with $t = 2m-4,2m-5,\dots,m-1$ we deduce that $S^{m-1}\CU^*(t+1-m) \in \D$
for these $t$ for even $m$. In both cases we see that all desired twists of $S^{m-1}\CU^*$ lie in $\D$.
\end{proof}

We will need also the following

\begin{lemma}\label{sko}
The bundles $\CS\otimes\CU^*(k)$ for $0\le k\le 2m-4$
as well as the bundles $\CS\otimes S^2\CU^*(k)$ for $0\le k\le 2m-5$ lie in the triangulated category
generated by the Lefschetz collection $\{S^l\CU^*(k)\ |\ (k,l)\in\Upso_{2m+1}\} \cup \{\CS(k)\ |\ 0\le k\le 2m-3\}$ in $\D^b(\OGr(2,W))$.
\end{lemma}
\begin{proof}
Consider the complex
$$
0 \to \CS(-1) \to \BS\otimes\CO_X \to \CS \to 0.
$$
By proposition~\ref{filtcs}
it is quasiisomorphic to $\CS\otimes\CU \cong \CS\otimes\CU^*(-1)$. Twisting it by $1,2,\dots,2m-3$
we deduce that $\CS\otimes\CU^*(k)$ lies in the desired triangulated category for $0\le k\le 2m-4$.
Further, consider the same complex tensored by $\CU^*(-1)$:
$$
0 \to \CS\otimes\CU^*(-2) \to \BS\otimes\CU^*(-1) \to \CS\otimes\CU^*(-1) \to 0.
$$
It is quasiisomorphic to
$\CS\otimes\CU\otimes\CU^*(-1) \cong \CS(-1) \oplus (\CS\otimes S^2\CU^*(-2))$.
Twisting it by $2,3,\dots,2m-3$ we deduce that $\CS\otimes S^2\CU^*(k)$ lies in the desired
triangulated category for $0\le k\le 2m-5$.
\end{proof}

The last preparatory result is the following:

\begin{lemma}\label{xphio}
For any two-dimensional subspace $\langle w_1,w_2\rangle \subset W \cong W^* = H^0(\OGr(2,W),\CU^*)$
such that the quadratic form $q$ is nondegenerate on $\langle w_1,w_2\rangle$, the zero locus of the corresponding section
$\phi = \phi_{w_1,w_2} \in H^0(\OGr(2,W),\CU^*\oplus\CU^*)$ on $X = \OGr(2,W)$
is the isotropic Grassmannian $X_{w_1,w_2} = \OGr(2,\lan w_1,w_2\ran^\perp) \subset \OGr(2,W) = X$.
Moreover, we have the following resolution of the structure sheaf $\CO_{X_\phi}$ on $X$:
\begin{equation}\label{resxphio}
0 \to \CO_X(-2) \to \CU(-1)\oplus \CU(-1) \to \CO_X(-1)^{\oplus 3} \oplus S^2\CU \to
\CU\oplus\CU \to \CO_X \to i_{\phi*}\CO_{X_{w_1,w_2}} \to 0,
\end{equation}
where $i_{\phi*}:X_{w_1,w_2} \to X$ is the embedding.
\end{lemma}
\begin{proof}
The first part is evident.
For the second part we note that any such section $\phi = \phi_{w_1,w_2}$ of $\CU^*\oplus\CU^*$
is regular since $\dim X_{w_1,w_2} = 2(n-2)-7 = \dim X - 4$, so the sheaf $i_{\phi*}\CO_{X_{w_1,w_2}}$
admits a Koszul resolution which takes form~(\ref{resxphis}).
\end{proof}

Now we are ready for the proof of the theorem.
We use induction in $m$. The base of induction, $m=2$, is clear.
Indeed, in this case $X = \OGr(2,W) = \PP^3$ and the Lefschetz collection
takes form $(\CO_{\PP^3},\CO_{\PP^3}(1),\CO_{\PP^3}(2),\CO_{\PP^3}(3))$
(since $\CO_X(1) = \CO_{\PP^3}(2)$ and $\CS = \CO_{\PP^3}(1)$)
which is well known to be full.

Now assume that the fullness of the corresponding Lefschetz collection is already proved for $m-1$.
Assume also that the Lefschetz collection for $m$ is not full. Then by~\cite{B} there exists an object $F\in\D^b(X)$,
right orthogonal to all bundles in the collection. Then by lemma~\ref{tsigmao} and lemma~\ref{sko} we have
$$
\begin{array}{rcl}
H^\bullet(X,S^l\CU(-k)\otimes F) = 0 &\qquad& \text{for all $(k,l)\in\TUpso_{2m-2}$}\\
H^\bullet(X,\CS^*(-k)\otimes F) = 0 &\qquad& \text{for all $0\le k\le 2m-3$}\\
H^\bullet(X,\CS^*\otimes\CU(-k)\otimes F) = 0 &\qquad& \text{for all $0\le k\le 2m-4$}\\
H^\bullet(X,\CS^*\otimes S^2\CU(-k)\otimes F) = 0 &\qquad& \text{for all $0\le k\le 2m-5$}
\end{array}
$$
Let us check that $i_\phi^*F = 0$ for any $\phi = \phi_{w_1,w_2}$ like in lemma~\ref{xphio}.
The same argument as in section~\ref{s_sgr} shows that $i_\phi^* F$ lies in the right orthogonal
to the subcategory of $\D^b(X_\phi)$ generated by the exceptional collection $\{S^l\CU^*(k)\ |\ (k,l)\in\Upso_{2m-1}\}$.
Similarly, tensor the resolution~(\ref{resxphio}) by $\CS^*(-k)\otimes F$ for $0\le k\le 2m-5$.
It follows that the cohomology on $X$ of the first five terms of the above complex vanishes.
Therefore we have
$$
\RHom_{X_\phi}(\CS(k),i_\phi^*F) = H^\bullet(X_\phi,\CS^*(-k)\otimes i_\phi^*(F)) = 0
\qquad\text{for $0\le k\le 2m-5$}
$$
Thus $i_\phi^* F$ lies in the right orthogonal to the subcategory of $\D^b(X_\phi)$
generated by the exceptional collection $\{\CS(k)\ |\ 0\le k\le 2m-5\}$.
But by the induction hypothesis the collection
$\{S^l\CU^*(k)\ |\ (k,l)\in\Upso_{2m-1}\} \cup \{\CS(k)\ |\ 0\le k\le 2m-5\}$
on $X_\phi$ is full, hence $i_\phi^*F = 0$.
So we conclude by the following

\begin{lemma}
If for $F\in\D^b(X)$ we have $i_\phi^*F = 0$ for any two-dimensional subspace
$\langle w_1,w_2\rangle \subset W$ such that the quadratic form $q$ is nondegenerate
on $\langle w_1,w_2\rangle$ then $F = 0$.
\end{lemma}
\begin{proof}
Assume that $F \ne 0$. Let $q$ be the maximal integer such that $\CH^q(F)\ne 0$,
take a point $x\in\supp\CH^q(F)$ and choose $\langle w_1,w_2\rangle \subset W$
such that $x\in X_\phi$ (this is equivalent to the orthogonality of $w_1$ and $w_2$
with the 2-dimensional subspace of $W$ corresponding to $x\in X = \OGr(2,W)$).
Since the functor $i_\phi^*$ is left-exact it easily follows that $\CH^q(i_\phi^* F) \ne 0$, so $i_\phi^*F \ne 0$.
\end{proof}

Thus we have proved that the desired collection is indeed full.

\begin{remark}
The same arguments allow to construct an exceptional collection on $\OGr(2,W)$
in case of even $n = \dim W$. Explicitly, if $n = 2m$, one can show that
the Lefschetz collection
\begin{equation*}\label{lcgro2m}
\left(
\begin{array}{rrrrrrr}
\CS_- & \CS_-(1) & \dots & \CS_-(m-1) & \CS_-(m-1) & \dots & \CS_-(2m-4) \\
\CS_+ & \CS_+(1) & \dots & \CS_+(m-1) & \CS_+(m-1) & \dots & \CS_+(2m-4) \\
S^{m-2}\CU^* & S^{m-2}\CU^*(1) & \dots & S^{m-2}\CU^*(m-2) \\
S^{m-3}\CU^* & S^{m-3}\CU^*(1) & \dots & S^{m-3}\CU^*(m-2) & S^{m-3}\CU^*(m-1) & \dots & S^{m-3}\CU^*(2m-4) \\
\vdots\quad & \vdots\quad & & \vdots\qquad\qquad & \vdots\qquad & & \vdots\qquad\qquad \\
\CU^* & \CU^*(1) & \dots & \CU^*(m-2) & \CU^*(m-1) & \dots & \CU^*(2m-4) \\
\CO & \CO(1) & \dots & \CO(m-2) & \CO(m-1) & \dots & \CO(2m-4)
\end{array}\right)
\end{equation*}
is exceptional. However, this collection is definitely not full
(its length is by 1 less then the rank of the Grothendieck group $K_0(\OGr(2,W))$).
So, to obtain a full collection one should add one more exceptional object
to the first block $(\CO,\CU^*,\dots,S^{m-2}\CU^*,\CS_+,\CS_-)$
of the above collection. This indeed can be done, but the object in question
turns out to be a complex (not a pure sheaf), so the picture became
more complicated. However, this is not quite satisfactory.
It seems that there should exist a better full exceptional collection.
\end{remark}


\begin{thebibliography}{XXXX}

\bibitem[Bei]{Bei} A. Beilinson,
{\em Coherent sheaves on $\PP^{n}$ and problems in linear algebra}, (Russian)
Funktsional. Anal. i Prilozhen. 12 (1978), no. 3, 68--69.,

\bibitem[B]{B} A. Bondal,
{\em Representations of associative algebras and coherent sheaves},
(Russian)  Izv. Akad. Nauk SSSR Ser. Mat. {\bf 53} (1989), no. 1, 25--44;
translation in  Math. USSR-Izv. {\bf 34}  (1990),  no. 1, 23--42.

\bibitem[Bou]{Bou} N. Bourbaki,
{\em \'El\'ements de math\'ematique. Premi\`ere partie: Les structures fondamentales de l'analyse.
Livre II: Alg\`ebre. Chapitre 9: Formes sesquilin\'eaires et formes quadratiques}, (French)
Actualit\'es Sci. Ind. no. 1272, Hermann, Paris, 1959, 211 pp. (1 insert).


\bibitem[D]{D} Demazure M.,
{\em A very simple proof of Bott's theorem},
Invent. Math. {\bf 33} (1976), no.~3, 271--272.

\bibitem[GR]{GR} A. Gorodentsev, A. Rudakov,
{\em Exceptional vector bundles on projective spaces},
Duke Math. J. {\bf 54} (1987), no.~1, 115--130.

\bibitem[Ka]{Ka} M. Kapranov,
{\em On the derived categories of coherent sheaves on some homogeneous spaces},
Invent. Math.  {\bf 92} (1988), no. 3, 479--508.

\bibitem[K1]{K1} A. Kuznetsov,
{\em Homological Projective Duality}, preprint math.AG/0507292.

\bibitem[K2]{K2} A. Kuznetsov,
{\em Lefschetz decompositions and Categorical resolutions of singularities},
preprint math.AG/0609240

\bibitem[K3]{K3} A. Kuznetsov,
{\em Homological projective duality for Grassmannians of lines},
preprint math.AG/0610957.

\bibitem[O]{O} D. Orlov,
{\em Projective bundles, monoidal transformations, and derived categories of coherent sheaves},
(Russian)  Izv. Ross. Akad. Nauk Ser. Mat. {\bf 56} (1992), no. 4, 852--862;
translation in  Russian Acad. Sci. Izv. Math. {\bf 41} (1993), no. 1, 133--141.

\bibitem[S1]{S1} A. Samokhin,
{\em The derived category of coherent sheaves on $LG_3^C$},
(Russian) Uspekhi Mat. Nauk {\bf 56} (2001), no. 3(339), 177--178;
translation in Russian Math. Surveys {\bf 56} (2001), no. 3, 592--594.

\bibitem[S2]{S2} A. Samokhin,
{\em Some remarks on the derived categories of coherent sheaves on homogeneous spaces},
preprint.


\end{thebibliography}
\end{document}